\documentclass[11pt,reqno]{amsart}      
\usepackage[english]{babel}
\usepackage[utf8x]{inputenc}
\usepackage[T1]{fontenc}
\usepackage[letterpaper, margin=1.2in]{geometry}
\usepackage{mathtools}
\usepackage[dvipsnames]{xcolor}
\usepackage{hyperref}
\usepackage{dsfont}
\usepackage{enumerate}
\usepackage{dsfont}
\usepackage[pdftex]{graphicx}
\usepackage{caption}
\usepackage{subcaption}
\usepackage{tikz}
\usetikzlibrary{positioning,decorations.pathreplacing,shapes.misc}
\tikzset{cross/.style={cross out, draw=black, fill=none, minimum size=2*(#1-\pgflinewidth), inner sep=0pt, outer sep=0pt}, cross/.default={2pt}}

\allowdisplaybreaks

\usepackage{amsmath, amsthm, amssymb}
\usepackage{latexsym,graphicx,multicol,mathrsfs,xypic}

\newtheorem{theorem}{Theorem}[section]
\newtheorem{definition}[theorem]{Definition}
\newtheorem{proposition}[theorem]{Proposition}
\newtheorem{lemma}[theorem]{Lemma}
\newtheorem{corollary}[theorem]{Corollary}

\newtheorem{remark}[theorem]{Remark}

\newtheorem*{convention}{Convention}
\newtheorem*{proposition*}{Proposition}

\DeclareMathOperator{\var}{var}
\DeclareMathOperator{\cov}{cov}

\DeclareMathOperator{\supp}{supp}

\usepackage{autonum}

\begin{document}

\title[Decay of correlations for the canonical ensemble]{Decay of correlations and uniqueness of the infinite-volume Gibbs measure of the canonical ensemble of 1d-lattice systems}

\author{Younghak Kwon}
\address{Department of Mathematics, University of California, Los Angeles}
\email{yhkwon@math.ucla.edu}

\author{Georg Menz}
\address{Department of Mathematics, University of California, Los Angeles}
\email{gmenz@math.ucla.edu}


\date{\today}

\begin{abstract}
We consider a one-dimensional lattice system of unbounded, real-valued spins with arbitrary strong, quadratic, finite-range interaction. We show the equivalence of correlations of the grand canonical (gce) and the canonical ensemble (ce). As a corollary we obtain that the correlations of the ce decay exponentially plus a volume correction term. Then, we use the decay of correlation to verify a conjecture that the infinite-volume Gibbs measure of the ce is unique on a one-dimensional lattice. For the equivalence of correlations, we modify a method that was recently used by the authors to show the equivalence of the ce and the gce on the level of thermodynamic functions. In this article we also show that the equivalence of the ce and the gce holds on the level of observables. One should be able to extend the methods and results to graphs with bounded degree as long as the gce has a sufficient strong decay of correlations. 
\end{abstract}

\maketitle

\section{Introduction}

The broader scope of this article is the study of phase transitions. A phase transition occurs if a microscopic change in a parameter leads to a fundamental change in one or more properties of the underlying physical system. The most well-known phase transition is when water becomes ice. Many physical systems, non-physical systems and mathematical models have phase transitions. For example, liquid to gas phase transitions are known as vaporization. Solid to liquid phase transitions are known as melting. Solid to gas phase transitions are known as sublimation. More examples are the phase transition in the 2-d Ising model (see for example~\cite{Sel16}), the Erd\"os-Renyi phase transition in random graphs (see for example~\cite{ErRe60},~\cite{ErRe61} or~\cite{KrSu13}) or phase transitions in social networks (see for example~\cite{Fro07}). \\

In this article, we study a one-dimensional lattice system of unbounded real-valued spins. The system consists of a finite number of sites~$i \in \Lambda \subset \mathbb{Z}$ on the lattice~$\mathbb{Z}$. For convenience, we assume that the set~$\Lambda$ is given by~$\left\{1, \ldots, N \right\}$. At each site~$i \in \Lambda$ there is a spin~$x_i$. In the Ising model the spins can take on the value~$0$ or~$1$. In this study, spins~$x_i \in \mathbb{R}$ are real-valued and unbounded. A configuration of the lattice system is given by a vector~$x \in \mathbb{R}^{K}$. The energy of a configuration~$x$ is given by the Hamiltonian~$H: \mathbb{R}^K \to \mathbb{R}$ of the system. For the detailed definition of the Hamiltonian~$H$ we refer to Section~\ref{s_setting_and_main_results}. We consider arbitrary strong, quadratic, finite-range interaction.\\

We consider two different ensembles: The first ensemble is the grand canonical ensemble (gce) which is given by the finite-volume Gibbs measures
\begin{align}
  \mu^\sigma(dx) = \frac{1}{Z} \exp\left( \sigma \sum_{i=1}^N x_i- H(x) \right) dx.
\end{align}
Here,~$Z$ is a generic normalization constant making the measure~$\mu^\sigma$ a probability measure. The constant~$\sigma \in \mathbb{R}$ is interpreted as an external field. The second ensemble is the canonical ensemble (ce). It emerges from the gce by conditioning on the mean spin (see Remark~\ref{r_relation_gce_ce})
\begin{align}
  m = \frac{1}{N} \sum_{i =1}^N x_i.
\end{align}
The ce is given by the probability measure
\begin{align}
  \mu_m(dx)& = \frac{1}{Z} \  \mathds{1}_{\left\{ \frac{1}{N}\sum_{i=1}^N x_i= m \right\}  } (x)\ \exp( - H(x)) \mathcal{L}^{N-1}(dx),
\end{align} 
where~$\mathcal{L}^{N-1}$ denotes the~$(N-1)$-dimensional Hausdorff measure.\\

There are many different ways to characterize phase transitions. In this article we use the convention that an ensemble has no phase transition if the associated infinite-volume Gibbs measure is unique. On the two-dimensional lattice the gce has a phase transition (see for example~\cite{Pei36}). This is not the case if one considers the gce on the one-dimensional lattice. There, the gce does not have a phase transition if the interaction is finite-range or decays fast enough (see for example~\cite{Isi25,Dob68,Dob74,Rue68,MeNi14}). It is a natural question if those results extend from the gce to the ce. This is a non-trivial question since there are known examples where the gce has no phase transition but the ce has (see for example~\cite{ScSh96,BiChKo02,BiChKo03}). For example, in~\cite{BiChKo03} the following situation is considered for the 2-d Ising model (cf.~Introduction, Theorem 1.1, Theorem 1.2 and Corollary~1.3 in~\cite{BiChKo03}). If the inverse temperature~$\beta$ is above the critical inverse temperature~$\beta_c$, it is known that in the grand canonical ensemble all the excess mass~$\Delta$ is concentrated in one large droplet determined by the so-called Wulff shape. This is not the case for the canonical ensemble which exhibits a phase transition: If the excess mass~$\Delta$ is smaller than a critical mass~$\Delta_c$ then all the mass is contained in droplets of logarithmic size. However, if the excess mass~$\Delta > \Delta_c$ almost all the mass is contained in a macroscopic droplet following the Wulff shape.  \\

The results of Cancrini, Martinelli and Roberto~\cite{CMR02} suggest that if the spins are discrete, i.e.~$\left\{0,1 \right\}$-valued, there is no phase transition for the ce on a one-dimensional lattice with nearest-neighbor interaction (see also the introduction of~\cite{KwMe18}). In this article we consider this question for real-valued and unbounded spins. Considering unbounded spins is harder because we lose compactness and one cannot transfer the arguments from the discrete case. It was conjectured in~\cite{KwMe18} that the infinite-volume Gibbs measure of the canonical ensemble is unique on the one-dimensional lattice. A first step toward verifying this conjecture was taken in~\cite{KwMe18}. There, it was shown that the gce and the ce are equivalent. This indicates that the gce and the ce should share similar properties. In this article, we show that the conjecture is indeed true: There is no phase transition for the ce on a one-dimensional lattice with finite-range interaction (see Theorem~\ref{p_uniqueness_iv_gibbs_ce} from below).\\

For the proof of Theorem~\ref{p_uniqueness_iv_gibbs_ce} we follow a standard argument (see e.g.~\cite{Yos03} or~\cite{Me13}) which is based on an ingredient of its own interest: decay of correlations. Decay of correlations is a classical ingredient when deducing the uniqueness of the infinite-volume Gibbs measure in lattice systems (see for example~\cite{DoSh87,MaBe99,MeNi14} for the gce). In many lattice systems decay of correlations is one of many equivalent mixing conditions, including the Dobrushin-Sloshman mixing condition (see e.g.~\cite{DoSh85,DoSh87,DoWa90,MaBe99,Yos03,HeMe16}). It is known that for one-dimensional lattice systems the correlations of the gce decay exponentially fast (see for example~\cite{Zeg96,MeNi14}, references therein and Theorem~\ref{p_decay_of_correlations_gce} below). In this article, we extend results of~\cite{MR00} for discrete-spin systems to unbounded continuous-spins (see also Section 3.2 in~\cite{CMR02}). More precisely, we show that the correlations of the ce and gce are equivalent (see Theorem~\ref{p_equivalence_decay_of_correlations} below). As a direct consequence we get that the correlations of the ce decay exponentially plus a volume correction term (see Theorem~\ref{p_exponential_decay_of_correlations_of_canonical_ensemble} below). We expect that this is a manifestation of a more general principle for systems on general lattices: If the gce has sufficient decay of correlations then the equivalence of correlations of the gce and ce should hold (see also~\cite{CaNiTr04}). In~\cite{KwMe18b}, the decay of correlations of the ce (cf.~Theorem~\ref{p_exponential_decay_of_correlations_of_canonical_ensemble}) is an important ingredient in the proof that the ce with arbitrary strong, ferromagnetic, finite-range interaction satisfies a uniform log-Sobolev inequality on the one-dimensional lattice.\\

For the proof of Theorem~\ref{p_exponential_decay_of_correlations_of_canonical_ensemble} we apply a similar strategy as for the proof of the equivalence of ensembles in~\cite{KwMe18}. We show that the decay of correlations can be transferred from the gce to the ce. The argument is robust and should apply to more general situations. The proof does not use that the lattice is one-dimensional. Instead, it only uses that the grand canonical ensemble on a one-dimensional lattice has an uniform exponential decay of correlations (see for example~\cite{MeNi14} and~\cite{Zeg96}). Under the assumption of fast enough decay of correlation, one should be able to use similar calculations to deduce decay of correlations of the ce for spin systems on arbitrary graphs, as long as the degree is uniformly bounded and the interaction has finite range. However, we only consider the one-dimensional lattice with finite-range interaction because less notational burden is better for explaining ideas and presenting the calculations. \\

There are many different aspects of equivalence of ensembles. For further background, we refer the reader to~\cite{StZe91,LePfSu94,Ge95,Ad06,Tou15}. After deducing the main results of this article but before submitting the manuscript, the authors heard of another approach to the equivalence of ensembles which was simultaneously and independently developed in~\cite{CaOl17}. For the presentation of our main results let us follow the exposition of~\cite{CaOl17}. There, the equivalence of the ce and the micro-canonical ensembles was deduced for classical particle systems via a combination of an Edgeworth expansion, a local central limit theorem and a local large deviations principle. Equivalence of ensembles exists on the level of thermodynamic functions, on the level of observables and on the level of correlations. In~\cite[Theorem~2]{KwMe18}, we showed a statement from which the equivalence of ensembles of the ce and gce on the level of thermodynamic functions follows. In this article, we will show that the equivalence of ensembles of the ce and gce also holds on the level of observables (see Theorem~\ref{p_equivalence_of_the_free_energies} below) and on the level of correlations (see~Theorem \ref{p_equivalence_decay_of_correlations}). In order for correlations~$\cov(f,g)$ to be equivalent, the involved functions~$f$ and~$g$ have to be local. If the functions are not local, the correlations and fluctuations, i.e.~covariances and variances, depend on the ensemble. However, there is still a nice relation between the fluctuations expressed by the Lebowitz-Percus-Verlet formula~\cite{LePeVe67}. This formula was rigorously deduced for particle systems in~\cite{CaOl17}. It would be very interesting to deduce the Lebowitz-Percus-Verlet formula in our setting which are unbounded spin systems. As discussed in~\cite{CaOl17}, their general strategy seems to apply to our setting. However, only a general strategy is given without explicit results like e.g.~the decay of correlation for the canonical ensemble. \\

Let us mention some more open questions and problems:
\begin{itemize}
\item Instead of considering finite-range interaction, is it possible to deduce similar results for infinite-range, algebraically decaying interactions? More precisely, is it possible to extend the results of~\cite{MeNi14} from the gce to the ce? Is the same algebraic order of decay sufficient, i.e.~of the order~$2+ \varepsilon$, or does one need a higher order of decay? \\[-1ex]

\item Is it possible to consider more general Hamiltonians? For example, our argument is based on the fact that the single-site potentials are perturbed quadratic, especially when we use the results of~\cite{KwMe18}. One would like to have general super-quadratic potentials as was for example used in~\cite{MeOt13}. Also, it would be nice to consider more general interactions than quadratic or pairwise interaction.\\[-1ex]
	
\item Is it possible to generalize the results to vector-valued spin systems?.
\end{itemize}

We conclude the introduction by giving an overview over the article. In Section~\ref{s_setting_and_main_results}, we introduce the precise setting and present the main results. In Section~\ref{s_auxiliary_lemmas}, we provide several auxiliary results. In Section~\ref{s_proof_exp_decay}, we prove the equivalence of the the gce and the ce on the level of observables and correlations (cf.~Theorem~\ref{p_equivalence_observables} and Theorem~\ref{p_equivalence_decay_of_correlations}). We also prove the decay of correlations of the canonical ensemble (cf.~Theorem~\ref{p_exponential_decay_of_correlations_of_canonical_ensemble}) in Section~\ref{s_proof_exp_decay}. In Section~\ref{s_uniqueness_ivGm}, we show the uniqueness of the infinite-volume Gibbs measure (cf.~Theorem~\ref{p_uniqueness_iv_gibbs_ce}).

\section*{Conventions and Notation}

\begin{itemize}
\item The symbol~$T_{(k)}$ denotes the term that is given by the line~$(k)$.
\item We denote with~$0<C<\infty$ a generic uniform constant. This means that the actual value of~$C$ might change from line to line or even within a line.
\item Uniform means that a statement holds uniformly in the system size~$|\Lambda|$, the mean spin~$m$, the boundary~$x^{\mathbb{Z}\backslash \Lambda}$ and the external field~$s$.
\item $a \lesssim b$ denotes that there is a uniform constant~$C$ such that~$a \leq C b$.
\item $a \sim b$ means that~$a \lesssim b$ and~$b \lesssim a$.
\item $\mathcal{L}^{k}$ denotes the $k$-dimensional Hausdorff measure. If there is no cause of confusion we write~$\mathcal{L}$.
\item $Z$ is a generic normalization constant. It denotes the partition function of a measure.  
\item For each~$N \in \mathbb{N}$,~$[N]$ denotes the set~$\left\{ 1, \ldots N \right\}$.
\item For a vector~$x \in \mathbb{R}^{\mathbb{Z}}$ and a set~$A \subset \mathbb{Z}$,~$x^A \in \mathbb{R}^{A}$ denotes the vector $ (x^A)_{i} = x_i$ for all~$i \in A$.
\item For a function~$f : \mathbb{R}^{\mathbb{Z}} \to \mathbb{C}$, denote~$\supp f $ by the minimal subset of~$\mathbb{Z}$ with $f(x) = f\left(x^{\supp f} \right)$. 
\item A function~$f : \mathbb{R}^{\mathbb{Z}} \to \mathbb{C}$ is said to be local if~$\supp f$ is finite.
\end{itemize}

\section{Setting and main results}
\label{s_setting_and_main_results}

We consider a system of unbounded continuous spins on the lattice~$\mathbb{Z}$. The formal Hamiltonian $H:\mathbb{R}^{\mathbb{Z}} \to \mathbb{R}$ of the system is defined as
\begin{align}\label{e_d_hamiltonian}
H(x) &= \sum_{i \in \mathbb{Z} } \left( \psi (x_i) + s_i x_i +\frac{1}{2}\sum_{j : \ 1 \leq |j-i| \leq R } M_{ij}x_i x_j \right) \\
& = \sum_{i \in \mathbb{Z}} \left( \psi_b (x_i) + s_ix_i + \frac{1}{2} \sum_{j : |j-i| \leq R} M_{ij}x_i x_j  \right),
\end{align}
where~$\psi(z) : = \frac{1}{2}z^2 + \psi_b (z)$ and~$M_{ii} : = 1$.
We assume the following:
\begin{itemize}
\item The function~$\psi_b: \mathbb{R} \to \mathbb{R}$ satisfies
 \begin{align}\label{e_nonconvexity_bounds_on_perturbation}
 |\psi_b|_{\infty} + |\psi'_b|_{\infty}  + |\psi''_b|_{\infty} < \infty. 
 \end{align}
It is best to imagine~$\psi(z) = \frac{1}{2}z^2 + \psi_b(z)$ as a double-well potential (see Figure~\ref{f_double_well}).  
\begin{figure}[t]
\centering
\begin{tikzpicture}
      \draw[->] (-2.5,0) -- (2.5,0) node[right] {$x$};
      \draw[->] (0,-1) -- (0,2.8) node[above] {};
      \draw[scale=1,domain=-1.9:1.9,smooth,variable=\x,blue] plot ({\x},{.3*\x*\x*\x*\x-.7*\x*\x-.3*\x+.2});
      
\node[align=center, below, blue] at (1.9, 1.9) {$\psi(x)$};      
\end{tikzpicture}
\caption{Example of a single-site potential~$\psi$}\label{f_double_well}
\end{figure}
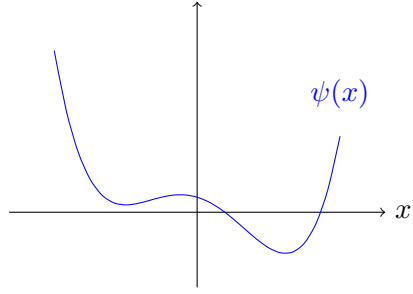

\item The interaction is symmetric i.e. 
\begin{align}
M_{ij} = M_{ji}  \qquad &\text{for all distinct} \ i, j \in \mathbb{Z}. 
\end{align}

\item The fixed, finite number~$R \in \mathbb{N}$ models the range of interactions between the particles in the system i.e.~it holds that~$M_{ij}=0$ for all~$i,j$ such that~$|i-j|> R$. \\

\item The matrix~$M=(M_{ij})$ is strictly diagonal dominant i.e. for some~$\delta>0$, it holds for any~$i \in \mathbb{Z}$ that
\begin{align}
\sum_{ 1 \leq |j-i| \leq R} |M_{ij}| + \delta \leq M_{ii} = 1. 
\end{align}

\item The vector~$s = (s_i) \in \mathbb{R}^{\mathbb{Z}}$ is arbitrary. It models the interaction with an inhomogeneous external field. Because the interaction is quadratic, this term also models the interaction of the system with the boundary. 
\end{itemize}



Let us consider a finite sublattice~$\Lambda \subset \mathbb{Z}$. Given boundary values ~$x^{\mathbb{Z} \backslash \Lambda} \in \mathbb{R}^{\mathbb{Z}\backslash \Lambda}$ we define the finite volume Hamiltonian~$H: \mathbb{R}^{\Lambda} \to \mathbb{R}$ as (using a small abuse of notation) 
\begin{align} 
  H(x^{\Lambda}) &: =H(x^{\Lambda}, x^{{\mathbb{Z}\backslash \Lambda}}) \\
  &= \sum_{i \in \Lambda } \Big( \psi (x_i) +  s_i  x_i +\frac{1}{2}\sum_{j : \ 1 \leq |j-i| \leq R } M_{ij}x_i x_j \Big) \\
  & = \sum_{i \in \Lambda } \left( \psi (x_i) +  \left( s_i + \frac{1}{2}\sum_{j\in \mathbb{Z} \backslash \Lambda : \ 1 \leq |j-i| \leq R } M_{ij} x_j  \right) x_i +\frac{1}{2}\sum_{j\in \Lambda : \ 1 \leq |j-i| \leq R } M_{ij}x_i x_j \right) . \label{e_def_finite_volume_Hamiltonian}
\end{align}
We want to point out that in~\eqref{e_def_finite_volume_Hamiltonian} the boundary values~$x^{\mathbb{Z}\backslash \Lambda}$ just modify the external field that is seen by a particular spin~$x_i$.

\begin{definition}
The gce~$\mu^{\Lambda, \sigma}$ associated to the Hamiltonian~$H$ with boundary values~$x^{\mathbb{Z} \backslash \Lambda}$ is the probability measure on~$\mathbb{R}^{\Lambda}$ given by the Lebesgue density
\begin{align} \label{d_gc_ensemble}
\mu^{\Lambda, \sigma} \left(dx^{\Lambda}\right) : = \frac{1}{Z} \exp\left( \sigma \sum_{k \in \Lambda} x_k  - H(x^{\Lambda}, x^{\mathbb{Z} \backslash \Lambda} ) \right) dx^{\Lambda}, 
\end{align}
where~$dx^{\Lambda}$ denotes the Lebesgue measure on~$\mathbb{R}^{\Lambda}$. The ce~$\mu_m ^{\Lambda}$ is the probability measure on
\begin{align}\label{e_d_X_lambda_m}
X_{\Lambda,m} : = \left \{ x^{\Lambda} \in \mathbb{R}^{\Lambda} : \ \frac{1}{|\Lambda|} \sum_{k \in \Lambda} x_k =m  \right \} \subset \mathbb{R}^{\Lambda}
\end{align}
with density
\begin{align} \label{d_ce}
\mu_m ^{\Lambda} (dx^{\Lambda}) : = \frac{1}{Z} \mathds{1}_{ \left\{ \frac{1}{|\Lambda|} \sum_{k \in \Lambda} x_k =m \right\}}\left(x^{\Lambda}\right) \exp\left( - H(x^{\Lambda}, x^{\mathbb{Z} \backslash \Lambda} ) \right) \mathcal{L}^{|\Lambda|-1}(dx^{\Lambda}). \label{d_ce}
\end{align}
Here~$\mathcal{L}^{|\Lambda|-1}(dx)$ denotes the~$(|\Lambda|-1)$-dimensional Hausdorff measure supported on~$X_{\Lambda, m}$. 
\end{definition}

\begin{remark} \label{r_relation_gce_ce} The ce~$\mu_m ^{\Lambda}$ emerges from the gce~$\mu^{\Lambda,\sigma}$ by conditioning on the mean spin
\begin{align} \label{e_spin_restriction}
\frac{1}{|\Lambda|} \sum_{ k \in \Lambda} x_k = m.
\end{align}
More precisely, given~\eqref{e_spin_restriction}, the term~$\sigma \sum_{k \in \Lambda} x_k$ inside the exponential in~\eqref{d_gc_ensemble} acts like a constant and hence is cancelled out with the normalization constant~$Z$ as follows:
\begin{align}
&\mu^{\Lambda, \sigma}  \left(dx^{\Lambda}  \mid \frac{1}{|\Lambda|}\sum_{k \in \Lambda} x_k =m \right) \\
&\qquad = \frac{1}{Z} \mathds{1}_{ \left\{ \frac{1}{|\Lambda|} \sum_{k \in \Lambda} x_k =m \right\}}\left(x^{\Lambda}\right) \exp\left(\sigma m |\Lambda| - H(x^{\Lambda}, x^{\mathbb{Z} \backslash \Lambda} ) \right) \mathcal{L}^{|\Lambda|-1}(dx^{\Lambda})\\
&\qquad = \frac{1}{\widetilde{Z}} \mathds{1}_{ \left\{ \frac{1}{|\Lambda|} \sum_{k \in \Lambda} x_k =m \right\}}\left(x^{\Lambda}\right) \exp\left( - H(x^{\Lambda}, x^{\mathbb{Z} \backslash \Lambda} ) \right) \mathcal{L}^{|\Lambda|-1}(dx^{\Lambda})\\
& \qquad = \mu_m ^{\Lambda} (dx^{\Lambda}).
\end{align}
Note that the ce~$\mu_m ^{\Lambda}$ does not depend on~$\sigma$, even though it emerged from the gce~$\mu^{\Lambda, \sigma}$.
\end{remark}
\medskip

To relate the external field~$\sigma$ of~$\mu^{\Lambda, \sigma}$ and the mean spin~$m$ of~$\mu_m ^{\Lambda}$ we further introduce following definitions.

\begin{definition}[The free energy of the gce] The free energy~$A_{gce} : \mathbb{R} \to \mathbb{R}$ of the gce~$\mu^{\Lambda, \sigma}$ is defined as
\begin{align}
A_{gce} (\sigma) : = \frac{1}{|\Lambda|}\ln \int_{\mathbb{R}^{\Lambda}} \exp \left( \sigma \sum_{k \in \Lambda} x_k - H\left( x^{\Lambda}, x^{\mathbb{Z} \backslash \Lambda}\right) \right)dx^{\Lambda}. \label{e_def_a_gce}
\end{align}
\end{definition}
\medskip

Let~$X= (X_k)_{k \in \Lambda}$ be a random variable distributed according to the gce~$\mu^{\Lambda,\sigma}$. A direct calculation yields
\begin{align}
\frac{d}{d\sigma}A_{gce} (\sigma) &= \frac{1}{|\Lambda|} \mathbb{E}_{\mu^{\Lambda, \sigma}} \left[\sum_{k \in \Lambda} X_k \right], \label{e_first_der_agce} \\
\frac{d^2}{d\sigma^2} A_{gce} (\sigma) &= \frac{1}{|\Lambda|} \var_{\mu^{\Lambda, \sigma}} \left(  \sum_{k \in \Lambda} X_k \right). \label{e_second_der_agce}
\end{align}
In~\cite{KwMe18}, the authors assumed
\begin{align} \label{e_needless}
\frac{1}{|\Lambda|} \var_{\mu^{\Lambda, \sigma}} \left( \sum_{k \in \Lambda} X_k \right) \gtrsim 1
\end{align}
to prove the strict convexity of~$A_{gce}$ in the sense that there is a constant~$C \in (0, \infty)$ with
\begin{align}
\frac{1}{C} \leq \frac{d^2}{d\sigma^2} A_{gce} \leq C.
\end{align}
There, the assumption~\eqref{e_needless} was attained in the case of attractive nearest-neighbor interactions (see~\cite[Lemma 1]{KwMe18}). In~\cite{KwMe18b},~\eqref{e_needless} was deduced again for the general finite-range interactions.

\begin{lemma}[Lemma 3.1 in~\cite{KwMe18b}] \label{l_variance_estimate} Let~$X= (X_k)_{k \in \Lambda}$ be a random variable distributed according to~$\mu^{\Lambda, \sigma}$. There is a constant~$C \in (0, \infty)$ which is uniform in the system size~$|\Lambda|$, the external fields~$s=(s_i)$ and~$\sigma$, and the boundary~$x^{\mathbb{Z} \setminus \Lambda}$ such that
\begin{align}
\frac{1}{C} \leq \frac{1}{|\Lambda|} \var_{\mu^{\Lambda, \sigma}} \left( \sum_{k \in \Lambda} X_k  \right)\leq C.
\end{align}
\end{lemma}
\medskip

A combination of~\eqref{e_second_der_agce} and Lemma~\ref{l_variance_estimate} yields the following corollary:

\begin{corollary}\label{c_strict_convexity_free_energy} The free energy~$A_{gce}$ of the gce~$\mu^{\Lambda, \sigma}$ is strictly convex in the sense that there is a constant~$C \in (0, \infty)$ independent of the system size~$|\Lambda|$, the external fields~$s=(s_i)$ and~$\sigma$, and the boundary~$x^{\mathbb{Z} \setminus \Lambda}$ such that
\begin{align}
\frac{1}{C} \leq \frac{d^2}{d\sigma^2} A_{gce} \leq C.
\end{align}
\end{corollary}
\medskip

Using the strict convexity of the free energy~$A_{gce}$, we relate the external field~$\sigma$ of~$\mu^{\Lambda, \sigma}$ and the mean spin~$m$ of~$\mu_m ^{\Lambda}$ as follows:

\begin{definition} \label{a_relation_sigma_m}
For each~$m \in \mathbb{R}$, we choose~$\sigma = \sigma(m) \in \mathbb{R}$ such that
\begin{align} \label{e_relation_sigma_m}
\frac{d}{d\sigma} A_{gce} (\sigma) = m,
\end{align}
or vice versa. Setting~$m_k := \mathbb{E}_{\mu^{\Lambda, \sigma}} \left[ X_k \right] $, this is equivalent to (see~\eqref{e_first_der_agce})
\begin{align} \label{e_m=sum_mi}
\frac{1}{|\Lambda|} \sum_{k \in \Lambda} m_k   = m.
\end{align}
\end{definition}
\medskip

Next, we define the free energy of the ce~$\mu_m^{\Lambda}$.

\begin{definition} [The free energy of the ce]
The free energy~$A_{ce} : \mathbb{R} \to \mathbb{R}$ of the ce~$\mu^{\Lambda}_m$ is defined (with slight abuse of notation) as 
\begin{align}
A_{ce}(\sigma) : = \frac{1}{|\Lambda|} \ln \int_{X_{\Lambda, m}} \exp \left( \sigma \sum_{k \in \Lambda} x_k - H\left( x^{\Lambda}, x^{\mathbb{Z} \backslash \Lambda} \right)  \right)\mathcal{L}^{|\Lambda|-1}(dx^{\Lambda}). \label{e_def_a_ce}
\end{align}
Because~$m = m(\sigma)$ can be viewed as a function of~$\sigma$, we simply write~$A_{ce} = A_{ce}(\sigma)$, instead of~$A_{ce} = A_{ce}(m , \sigma)$.
\end{definition}
\medskip

In this article we consider three different notions of equivalence of ensembles. The first one is the equivalence of the free energies in~$C^2$ and was deduced in~\cite{KwMe18}. 

\begin{theorem}[Theorem 2 in~\cite{KwMe18}] \label{p_equivalence_of_the_free_energies}
It holds that
\begin{align}
\lim_{\Lambda \to \mathbb{Z}} \left|A_{gce} - A_{ce} \right|_{C^2}=0,
\end{align}
where the convergence is uniform. More precisely, for each~$\varepsilon>0$, there is an integer~$N_0 \in \mathbb{N}$ independent of the external field~$s=(s_i)$ and the boundary~$x^{\mathbb{Z} \setminus \Lambda}$ such that for all~$\Lambda \subset{\mathbb{Z}}$ with~$|\Lambda| \geq N_0$,
\begin{align}
\sup_{\sigma \in \mathbb{R}} \left| A_{gce} (\sigma) -A_{ce} (\sigma) \right| &\lesssim \frac{1}{|\Lambda|},\\
\sup_{\sigma \in \mathbb{R}} \left|\frac{d}{d\sigma} A_{gce}(\sigma) - \frac{d}{d\sigma} A_{ce}(\sigma) \right| &\lesssim \frac{1}{|\Lambda|^{1-\varepsilon}}, \\
\sup_{\sigma \in \mathbb{R}} \left| \frac{d^2}{d\sigma^2} A_{gce}(\sigma) - \frac{d^2}{d\sigma^2} A_{ce}(\sigma) \right| &\lesssim \frac{1}{|\Lambda|^{\frac{1}{2}-\varepsilon}}.
\end{align}
\end{theorem}
\medskip

The second one is the equivalence of observables.

\begin{theorem}\label{p_equivalence_observables}
For a function~$f : \mathbb{R}^{\mathbb{Z}} \to \mathbb{C}$, denote~$\supp f $ by the minimal subset of~$\mathbb{Z}$ with $f(x) = f\left(x^{\supp f} \right)$. Let~$f: \mathbb{R}^{\mathbb{Z}} \to \mathbb{R}$ be a local function, i.e. a function with finite support. Then for each~$\varepsilon>0$, there exist constants~$N_0 \in \mathbb{N}$ and~$\tilde{C}=\tilde{C}(\varepsilon) \in (0, \infty)$ independent of the external fields~$s=(s_i)$, the boundary~$x^{\mathbb{Z} \backslash \Lambda}$, and the mean spin~$m$ such that for all~$\Lambda \subset \mathbb{Z}$ with~$|\Lambda| \geq N_0$ and~$\Lambda \supset \supp f$, it holds that
\begin{align}
\left| \mathbb{E}_{\mu^{\Lambda, \sigma}}\left[ f\right] - \mathbb{E}_{\mu^{\Lambda}_m} \left[ f \right] \right| \leq \tilde{C} \frac{\| \nabla f \|_{L^2 (\mu^{\Lambda, \sigma})} |\supp f|^{\frac{1}{2}}}{N^{ \frac{1}{2} - \varepsilon}}. \label{e_proof_equivalence_observables}
\end{align}
\end{theorem}
\medskip

\begin{remark}
Recall that for each~$m \in \mathbb{R}$, we choose~$\sigma= \sigma(m)$ by~\eqref{e_relation_sigma_m} or vice versa. With this definition, the statement of Theorem~\ref{p_equivalence_observables} can also be understood in the following sense: For each~$\varepsilon >0$, there exist constants~$N_0 \in \mathbb{N}$ and~$\tilde{C} = \tilde{C}(\varepsilon) \in (0, \infty)$ independent of the external fields~$s= (s_i)$ and~$\sigma$, and the boundary~$x^{\mathbb{Z} \setminus \Lambda}$ such that for all~$\Lambda \subset \mathbb{Z}$ with~$|\Lambda| \geq N_0$ and~$\Lambda \supset \supp f$,~\eqref{e_proof_equivalence_observables} holds.
\end{remark}
\medskip 
    
We give the proof of Theorem~\ref{p_equivalence_observables} in Section~\ref{s_proof_exp_decay}. The main idea is to relate the left hand side of~\eqref{e_proof_equivalence_observables} to the derivatives of free energy of suitable modified ensembles (see Lemma~\ref{l_second_derivative_covariance}). This trick is well known in the community, for example, the arguments presented in~\cite{MR00} rely on the same trick. \\

In this article, we will study decay of correlations of the ce~$\mu_m^{\Lambda}$. For that purpose let us, before we proceed, recall that for one-dimensional lattice systems the correlations of the gce decay exponentially fast (\cite[Lemma 6]{KwMe18}. See also~\cite[Theorem 1.4]{MeNi14}). 

\begin{theorem}[Lemma 6 in~\cite{KwMe18}] \label{p_decay_of_correlations_gce}
Let~$f, g : \mathbb{R}^{\mathbb{Z}} \to \mathbb{R}$ be local functions supported on~$\Lambda$. Then 
\begin{align}
&\left| \cov_{\mu^{\Lambda, \sigma}} \left( f, g \right) \right| \lesssim \|\nabla f\|_{L^2 (\mu^{\Lambda, \sigma})}\|\nabla g\|_{L^2 (\mu^{\Lambda, \sigma})} \exp \left( -C \text{dist} \left( \supp f , \supp g \right) \right). \label{e_decay_of_correlation_ce}
\end{align}
\end{theorem}
\medskip

The last notion of equivalence of ensembles is the equivalence of correlations.

\begin{theorem}\label{p_equivalence_decay_of_correlations}
Let~$f, g : \mathbb{R}^{\mathbb{Z}} \to \mathbb{R}$ be local functions. Denote~$C(f,g)$ by
\begin{align}
C (f,g) : &= \left\| \nabla \left( \left(f(X) - \mathbb{E}_{\mu^{\Lambda, \sigma}}\left[ f(X)\right]\right)\left(g(X) - \mathbb{E}_{\mu^{\Lambda, \sigma}}\left[ g(X)\right]\right)   \right)\right\|_{L^{4}(\mu^{\Lambda, \sigma})} \\
& \quad + \| \nabla f \|_{L^{4}(\mu^{\Lambda, \sigma})}\| \nabla g \|_{L^{4}(\mu^{\Lambda, \sigma})}.    \label{e_def_c}
\end{align}
Then for each~$\varepsilon>0$, there exist constants~$N_0$ and~$\tilde{C}=\tilde{C}( \varepsilon) \in (0, \infty)$ independent of the~$s=(s_i)$, the boundary~$x^{\mathbb{Z}\backslash \Lambda}$, and the mean spin~$m$ such that for all~$\Lambda \subset \mathbb{Z}$ with~$|\Lambda| \geq N_0$ and~$ \Lambda \supset \supp f$, it holds that
\begin{align}
& \left|  \cov_{\mu_m^{\Lambda}}  (f,g) - \cov_{\mu^{\Lambda, \sigma}}  (f,g) \right| \\
&\qquad \leq \tilde{C} \ C(f,g) \left( \frac{ |\supp f | + |\supp g |  }{|\Lambda|^{1- \varepsilon}} + \exp\left(-C\text{dist}\left( \supp f, \supp g \right) \right) \right)  . \label{e_equiv_dec_cor}
\end{align}
\end{theorem}
\medskip 

We give the proof of Theorem~\ref{p_equivalence_decay_of_correlations} in Section~\ref{s_proof_exp_decay}. As in Theorem~\ref{p_equivalence_observables}, the first step is to relate the left hand side of~\eqref{e_equiv_dec_cor} to the derivatives of free energy of suitable modified ensembles (see Lemma~\ref{l_second_derivative_covariance}). \\

A combination of Theorem~\ref{p_decay_of_correlations_gce} and Theorem~\ref{p_equivalence_decay_of_correlations} yields another main result of this article:

\begin{theorem}[Decay of correlations of the ce]\label{p_exponential_decay_of_correlations_of_canonical_ensemble}
Under the same assumptions as in Theorem~\ref{p_equivalence_decay_of_correlations}, it holds that
\begin{align} \label{e_decay_correlations_ce}
\left|\cov_{\mu_m ^{\Lambda}} \left( f, g \right) \right| \leq \tilde{C} \ C(f,g) \left( \frac{ |\supp f | + |\supp g |  }{|\Lambda|^{1- \varepsilon}} + \exp\left(-C\text{dist}\left( \supp f, \supp g \right) \right) \right).
\end{align}
\end{theorem}
\medskip

We give the proof of Theorem~\ref{p_exponential_decay_of_correlations_of_canonical_ensemble} in Section~\ref{s_proof_exp_decay}. 

\begin{remark}\label{p_compare_remark}
In~\cite{TsDi15} a similar result for classical particle systems was deduced. However, there are subtle and important differences between \cite[Theorem 2.5]{TsDi15} and Theorem~\ref{p_equivalence_decay_of_correlations}. Obviously, one statement is about lattice systems and the other one is about particle systems. The next difference is that Theorem~\ref{p_equivalence_decay_of_correlations} yields the decay of correlations for arbitrary local functions~$f$ and~$g$ whereas~\cite[Theorem 2.5]{TsDi15} only shows the decay of correlations for the two-point function, i.e.~setting~$f(x)=x_i$ and~$g(x)=x_j$. Another subtle difference is that for~\cite[Theorem 2.5]{TsDi15} one needs the validity of the cluster expansion for the ce. For deducing Theorem~\ref{p_equivalence_decay_of_correlations} on general lattices one would only need the decay of correlations for the \emph{gce}. Cluster expansions and decay of correlations are closely related. For example, both hold in the high temperature regime. However, studying and deriving properties for the gce is a lot easier than for the ce.
\end{remark}
\medskip 

\begin{remark} 
The statements of Theorem~\ref{p_equivalence_observables} and Theorem~\ref{p_exponential_decay_of_correlations_of_canonical_ensemble} can be understood as an extension of a classical result of Cancrini and Martinelli from discrete spins to unbounded continuous spins (see Theorem 4.1 and Section 7.3 in~\cite{MR00}). The bounds of Theorem~\ref{p_equivalence_observables} and Theorem~\ref{p_exponential_decay_of_correlations_of_canonical_ensemble} are sub-optimal compared to~\cite{MR00} in the following sense. The right hand side of Theorem~\ref{p_equivalence_observables} and of Theorem~\ref{p_exponential_decay_of_correlations_of_canonical_ensemble} scale like~$|\Lambda|^{\frac{1}{2}-\varepsilon}$ and~$|\Lambda|^{1- \varepsilon}$ respectively, while the estimates of~\cite{MR00} scale like~$|\Lambda|^{-1}$.
This sub-optimality is due to the unboundedness of the spins, which make the problem more difficult. There is also another important difference in the estimates. The bounds in~\cite{MR00} include~$L^{\infty}$ norms. In the setting of unbounded spins, such a bound would be of limited value. For example one cannot deduce decay of correlations of the spins, i.e. setting~$f(x) =x_i$ and~$g(x) = x_j$ (see Corollary~\ref{p_decay_two_point_ce}). Also, the derivation of one of our main results, namely Theorem~\ref{p_uniqueness_iv_gibbs_ce}, heavily relies on the decay of correlations of the spins. An interesting problem is if our bounds could be improved to get the same scaling in~$|\Lambda|$ as in~\cite{MR00}. The recent results of~\cite{CaOl17} for classical particle systems indicate that this is possible.
\end{remark}
\medskip 

In order to deduce decay of correlations via Theorem~\ref{p_exponential_decay_of_correlations_of_canonical_ensemble} one still needs to estimate the right hand side of~\eqref{e_decay_correlations_ce}. Fortunately, the right hand side only involves $L^k$ norms with respect to gce~$\mu^{\Lambda,\sigma}$ and not the ce~$\mu_m^{\Lambda}$. Given the fact that under sufficient decay of correlations the gce~$\mu^{\Lambda,\sigma}$ satisfies a uniform LSI and Poincar\'e inequality (cf.~\cite{HeMe16}), it is a lot easier to estimate $L^k$ norms with respect to the gce~$\mu^{\Lambda,\sigma}$ than with respect to the ce~$\mu_m^{\Lambda}$. For example, we have the following statement:

\begin{lemma}[Lemma 5 in~\cite{KwMe18}]\label{l_moment_estimate} 
For each~$i \in \Lambda$, define
\begin{align} \label{d_def_mi}
m_i : = \mathbb{E}_{\mu^{\Lambda, \sigma}} \left[ X_i \right].
\end{align}
Then for any~$k \geq 1$, there is a constant~$C(k)$ independent of the system size~$|\Lambda|$, the external fields~$s=(s_i)$, the boundary~$x^{\mathbb{Z} \backslash \Lambda}$, and the mean spin~$m$ such that
\begin{align}
\mathbb{E}_{\mu^{\Lambda, \sigma}} \left[ \left| X_i - m_i \right|^k \right] \leq C(k) \qquad \text{for all } i \in \Lambda.
\end{align}
\end{lemma}
\medskip 

Lemma~\ref{l_moment_estimate} was proved in~\cite{KwMe18} in the case of nearest-neighbor interactions. However, with a cosmetic change, the result can be generalized to the case of strong finite-range interactions. A combination of Theorem~\ref{p_exponential_decay_of_correlations_of_canonical_ensemble} and Lemma~\ref{l_moment_estimate} yields the following statement.
\begin{corollary}\label{p_decay_two_point_ce}
For given~$\varepsilon>0$, there exist constants~$N_0 \in \mathbb{N}$ and~$\tilde{C}=\tilde{C}(\varepsilon)$ independent of the external fields~$s=(s_i)$, the boundary~$x^{\mathbb{Z} \backslash \Lambda}$, and the mean spin~$m$ such that for any sublattice~$\Lambda \subset \mathbb{Z}$ with~$|\Lambda| \geq N_0$, it holds that for any~$i, j \in \Lambda$,
\begin{align}
\left| \cov_{\mu^{\Lambda}_m} \left( X_i, X_j \right) \right| \leq \tilde{C} \left( \frac{1}{|\Lambda|^{1-\varepsilon}} + \exp \left( - C |i-j| \right)\right). \label{e_can_exp_dec}
\end{align}
\end{corollary}
\medskip 

\begin{remark}
Comparing the estimate~\eqref{e_can_exp_dec} of Corollary~\ref{p_decay_two_point_ce} with the estimate~\eqref{e_decay_of_correlation_ce} of Theorem~\ref{p_decay_of_correlations_gce} one observes the additional term~$\frac{1}{|\Lambda|^{1-\varepsilon}}$ on the right hand side of~\eqref{e_can_exp_dec}. Eventually, this term could be improved to the order~$\frac{1}{|\Lambda|}$ but not further. For example, assuming the random variables~$X_1, \cdots, X_N$ are distributed according to a symmetric measure~$\nu$ with fixed mean spin~$m$, it holds that for all~$i,j \in [N]$ with~$i \neq j$ 
\begin{align}
\cov_{\nu} \left( X_1, X_2 \right) &= \cov_{\nu} \left( X_i, X_j \right).
\end{align}
Therefore we have
\begin{align}
\cov_{\nu} \left( X_1, X_2 \right) &= \frac{1}{N-1} \cov_{\nu} \left( X_1, X_2 + \cdots X_N \right) \\
&= \frac{1}{N-1} \cov_{\nu} \left( X_1, Nm - X_1 \right) = - \frac{1}{N-1} \var_{\nu} \left(X_1 \right).
\end{align}
\end{remark}
\medskip

Let us turn to the next main result of this article, namely the uniqueness of the infinite-volume Gibbs measure of the ce.
\begin{definition}[Infinte-volume Gibbs measure] Let~$\mu$ be a probability measure on~$\mathbb{R}^{\mathbb{Z}}$ with standard product Borel sigma-algebra. For any finite subset~$\Lambda \subset \mathbb{Z}$ we decompose the measure~$\mu$ into the conditional measure~$\mu(dx^{\Lambda}|x^{\mathbb{Z} \backslash \Lambda})$ and the marginal~$\bar \mu (dx^{\mathbb{Z} \backslash \Lambda})$. This means that for any test function~$f$ it holds
\begin{align}
 \int f(x) \mu (dx) = \int \int f(x^{\Lambda}, x^{\mathbb{Z}\backslash \Lambda}) \mu \left( dx^\Lambda| x^{\mathbb{Z} \backslash \Lambda} \right) \mu\left(dx^{\mathbb{Z} \backslash \Lambda}\right). 
\end{align}
We say that~$\mu$ is an infinite-volume Gibbs measure of the ce if the conditional measures~$\mu \left( dx^\Lambda| x^{\mathbb{Z} \backslash \Lambda} \right)$ are given by the finite volume ce~$\mu_m^{\Lambda}(dx^{\Lambda})$ given by Definition~\ref{d_ce} i.e.
\begin{align} \label{e_dlr}
\mu \left( dx^\Lambda| x^{\mathbb{Z} \backslash \Lambda} \right) = \mu_m^{\Lambda}(dx^{\Lambda}).
\end{align}
The equations of the last identity are called Dobrushin-Lanford-Ruelle (DLR) equations. 
\end{definition}
\medskip 

\begin{theorem}[Uniqueness of the infinite-volume Gibbs measure of the ce]\label{p_uniqueness_iv_gibbs_ce}
Let~$H$ be a Hamiltonian that satisfies the assumptions described at the beginning of this section. Then there is only one  infinite-volume Gibbs measure of the ce that satisfies the uniform bound
\begin{align}
\sup_{i \in \mathbb{Z}} \var_{\mu} (x_i)< \infty. \label{e_uniform_bound_uniqueness}
\end{align}
\end{theorem}
\medskip

We deduce Theorem~\ref{p_uniqueness_iv_gibbs_ce} in Section~\ref{s_uniqueness_ivGm}. The main ingredient in the proof is the decay of correlations (cf. Theorem~\ref{p_exponential_decay_of_correlations_of_canonical_ensemble}).

\begin{remark} In this article, we only show the uniqueness of infinite-volume Gibbs measure of the ce, not the existence. However, the authors of this article believe that with a cosmetic change, the existence should follow by a compactness argument (see for example~\cite{BeHoKr82}).
\end{remark}
\medskip

Let us comment on how the main results are connected. Theorem~\ref{p_equivalence_observables} is shown fast. Adapting the techniques used for proving Theorem~\ref{p_equivalence_observables}, Theorem~\ref{p_equivalence_decay_of_correlations} is also deduced. Theorem~\ref{p_exponential_decay_of_correlations_of_canonical_ensemble} is a consequence of Theorem~\ref{p_decay_of_correlations_gce} and Theorem~\ref{p_equivalence_decay_of_correlations}. Theorem~\ref{p_exponential_decay_of_correlations_of_canonical_ensemble} is used as a main ingredient for proving Theorem~\ref{p_uniqueness_iv_gibbs_ce}.

\section{Auxiliary Lemmas} \label{s_auxiliary_lemmas}

In this section we provide auxiliary results. Lemma~\ref{l_amgm} and Proposition~\ref{p_main_computation} were proved in~\cite{KwMe18} for lattice systems with nearest-neighbor interaction. However, it is not hard to see that the arguments in~\cite{KwMe18} can be generalized in a straight-forward manner to lattice systems with finite range interaction $R< \infty$, which is considered in this article. \\

Recall the definition~\eqref{d_gc_ensemble} of gce
\begin{align}
\mu^{\Lambda, \sigma}\left(dx^{\Lambda}\right) : = \frac{1}{Z} \exp\left( \sigma \sum_{k \in \Lambda} x_k -H\left(x^{\Lambda}, x^{\mathbb{Z} \backslash \Lambda}\right) \right) dx^{\Lambda}. 
\end{align}
The first statement is a direct consequence of Lemma~\ref{l_moment_estimate}.

\begin{lemma}[Lemma 5 in~\cite{KwMe18}]  \label{l_amgm}
For any finite set~$B_i \subset \Lambda$ and~$k \in \mathbb{N}$, it holds that
\begin{align}
\left|\mathbb{E}_{\mu^{\Lambda, \sigma}}\left[ \sum_{i_1 \in B_1 } \cdots \sum_{i_k \in B_k} \left( X_{i_1} - m_{i_1} \right) \cdots \left( X_{i_k} - m_{i_k} \right)  \right]  \right| \lesssim \left| B_1 \right| \cdots \left| B_k \right|.
\end{align}
\end{lemma}
\medskip

The next statement provides the fourth moment estimate.

\begin{lemma} \label{l_fourth_moment}
For any subset~$A$ of $\Lambda$, it holds that
\begin{align}
\mathbb{E}_{\mu^{\Lambda, \sigma}} \left[ \left( \sum_{i \in A} \left(X_i -m_i \right) \right)^4 \right] \lesssim |A|^2.
\end{align}
\end{lemma}
\medskip

\noindent \emph{Proof of Lemma~\ref{l_fourth_moment}.} \ We prove the case when~$\Lambda = A = \{1, \cdots, N\}$. General case follows from the similar argument. For each pair~$(i,j,k,l) \subset \{1, \cdots, N\}$ with~$i \leq j \leq k \leq l$ we have by Theorem~\ref{p_decay_of_correlations_gce} and Lemma~\ref{l_amgm} that
\begin{align}
&\left| \mathbb{E}_{\mu^{\Lambda, \sigma}} \left[ (X_i -m_i ) (X_j - m_j) (X_k -m_k) (X_l- m_l) \right] \right| \\
& \qquad = \left| \cov_{\mu^{\Lambda, \sigma}} \left( X_i, (X_j - m_j) (X_k -m_k) (X_l- m_l) \right)\right| \lesssim \exp \left( - C |i-j| \right).
\end{align}
Similarly, one gets
\begin{align}
&\left| \mathbb{E}_{\mu^{\Lambda, \sigma}} \left[ (X_i -m_i ) (X_j - m_j) (X_k -m_k) (X_l- m_l) \right] \right| \lesssim \exp \left( - C |k-l| \right),
\end{align}
and thus
\begin{align}
\left| \mathbb{E}_{\mu^{\Lambda, \sigma}} \left[ (X_i -m_i ) (X_j - m_j) (X_k -m_k) (X_l- m_l) \right] \right| \lesssim \exp \left( - C \max \left(|i-j|, |k-l| \right) \right).    
\end{align}
Therefore we conclude
\begin{align}
\mathbb{E}_{\mu^{\sigma}} \left[ \left( \sum_{i=1}^N Y_i \right)^4 \right] \lesssim \sum_{i \leq j \leq k \leq l} \mathbb{E}_{\mu^{\sigma}} \left[ Y_i Y_j Y_k Y_l \right] \lesssim \sum_{j \leq k} \sum_{d=0}^{N-1} 2(d+1) \exp \left( -Cd \right) \lesssim N^2.
\end{align}
\qed
\medskip

Lastly, let~$g$ be the density of the random variable
\begin{align}
\frac{1}{\sqrt{|\Lambda|}} \sum_{k \in \Lambda} \left( X_k - m \right) \overset{\eqref{e_m=sum_mi}}{=} \frac{1}{\sqrt{|\Lambda|}} \sum_{k \in \Lambda} \left( X_k - m_k \right),
\end{align}
where the random vector~$X=(X_k)_{k \in \Lambda}$ is distributed according to~$\mu^{\Lambda, \sigma}$.
The following proposition provides estimates for~$g(0)$.

\begin{proposition} [Proposition 1 in~\cite{KwMe18}] \label{p_main_computation}
For each~$\alpha>0$ and~$\beta > \frac{1}{2}$, there exist constants~$C \in (0, \infty)$ and~$N_0 \in \mathbb{N}$ independent of the external fields~$s=(s_i)$, the boundary~$x^{\mathbb{Z} \backslash \Lambda}$, and the mean spin~$m$ such that for all~$\Lambda$ with $|\Lambda| \geq N_0$, it holds that
\begin{align}
\frac{1}{C} \leq g(0) \leq C , \qquad 
\left| \frac{d}{d\sigma} g(0) \right| \lesssim |\Lambda|^{\alpha}  \qquad \text{and} \qquad 
\left| \frac{d^2}{d\sigma^2} g(0) \right| \lesssim |\Lambda|^{\beta}.
\end{align}
\end{proposition}

\section{Proof of Theorem~\ref{p_equivalence_observables}, Theorem~\ref{p_equivalence_decay_of_correlations} and Theorem~\ref{p_exponential_decay_of_correlations_of_canonical_ensemble}} \label{s_proof_exp_decay}

The main idea of the proof of Theorem~\ref{p_equivalence_observables} and Theorem~\ref{p_equivalence_decay_of_correlations} is to write the left hand sides of the statements in terms of~$L^k$ norms with respect to gce~$\mu^{\Lambda, \sigma}$ via Cram\'er's representation. Since the gce~$\mu^{\Lambda, \sigma}$ on a one-dimensional lattice satisfies a uniform LSI and Poincar\'e inequality (cf.~\cite{HeMe16}), the~$L^k$ norms with respect to the gce~$\mu^{\Lambda, \sigma}$ can be estimated relatively easily. \\

One difficulty of this argument is that we consider general local functions~$f$ and~$g$. Indeed, the case when~$f$ and~$g$ are point functions, i.e.~$f(x) =x_i$ and~$g(x)=x_j$ for some~$i, j \in \Lambda$ (cf. Corollary~\ref{p_decay_two_point_ce}), one can easily prove the theorem with the help of moment estimates given in Lemma~\ref{l_moment_estimate}. However, we overcome this difficulty by combining Lemma~\ref{l_moment_estimate} with Poincar\'e inequality and H\"{o}lder's inequality. For more details, we refer to Section~\ref{s_equiv_observables_proof_of_proposition} and Section~\ref{s_exp_decay_proof_of_proposition}. \\

Let us begin with a convention which will reduce our notational burden.

\begin{convention}
We assume that~$\Lambda= [N] = \left\{1, \cdots, N \right\}$. Moreover, if there is no source of confusion, we write (with some abuse of notations)~$\mu^{\sigma} : = \mu^{[N], \sigma}$, $\mu_m : = \mu^{[N]}_m$,~$x= x^{[N]}$ and~$H(x) = H\left(x^{[N]}, x^{\mathbb{Z} \backslash [N]}\right)$.
\end{convention}
\medskip

Let us introduce auxiliary notations and settings. Fix the external field~$\sigma$ and local functions~$f$,~$g$. Define modified gce~$\mu^{\sigma_i, \sigma_j} $ and ce~$\mu_m ^{\sigma_i, \sigma_j}$ depending on~$\sigma_i, \sigma_j \in \mathbb{R}$ by
\begin{align}
\mu^{\sigma_i, \sigma_j }  \left(dx\right) &: = \frac{1}{Z} \exp\left( \sigma \sum_{k=1}^N x_k + \sigma_i f\left(x\right) + \sigma_j g\left(x\right) -H(x) \right) dx, \label{e_modified_fg_gce}\\
\mu_m^{\sigma_i, \sigma_j }  \left(dx\right) &: = \frac{1}{Z} \mathds{1}_{ \left\{ \frac{1}{N} \sum_{k=1}^N x_k =m \right\}}(x) \exp\left( \sigma \sum_{k=1}^N x_k + \sigma_i f\left(x\right) + \sigma_j g\left(x\right) -H(x) \right) \mathcal{L}^{N-1}(dx). \label{e_modified_fg_ce}
\end{align}
Note in particular that we have
\begin{align} \label{e_relation_between_measures}
\mu^{0, 0}  = \mu ^{\sigma}, \quad \mu_m ^{0, 0} = \mu_m.
\end{align} 
The associated free energies are
\begin{align}
A_{gce}^{f,g}\left(\sigma_i, \sigma_j \right) &: = \frac{1}{N} \ln \int \exp\left( \sigma \sum_{k=1}^N x_k + \sigma_i f(x) + \sigma_j g(x) -H(x) \right) dx,  \label{d_agce_ij}  \\
A_{ce}^{f,g}\left(\sigma_i, \sigma_j \right)  &: = \frac{1}{N} \ln \int_{\left\{\frac{1}{N} \sum_{k=1}^N x_k =m  \right\}} \exp\left( \sigma \sum_{k=1}^N x_k + \sigma_i f(x) + \sigma_j g(x) -H(x) \right) \mathcal{L}^{N-1}(dx). \label{d_ace_ij}
\end{align}

There are two ways to interpret free energies. 
First, the following auxiliary lemma (Lemma~\ref{l_second_derivative_covariance}) connects the free energies~$A_{ce}^{f,g}$ and~$A_{gce}^{f,g}$ of the modified ensembles~$\mu^{\sigma_i, \sigma_j}$ and~$\mu_{m}^{\sigma_i, \sigma_j}$ with expectations and covariances with respect to the gce~$\mu^{\sigma}$ and the ce~$\mu_m$.

\begin{lemma} \label{l_second_derivative_covariance}
Let~$X$ and~$Y$ be real-valued random variables distributed according to the gce~$\mu^{\sigma}$ and the ce~$\mu_m$, respectively. Then it holds that
\begin{align}
&\left.\frac{d}{d\sigma_i } A_{gce}^{f,g} \right|_{ \sigma_i , \sigma_j =0} = \frac{1}{N} \mathbb{E}_{\mu^{\sigma}} \left[ f(X) \right], \label{e_agce_exp} \\
&\left.\frac{d}{d\sigma_i } A_{ce}^{f,g} \right|_{ \sigma_i , \sigma_j =0} = \frac{1}{N} \mathbb{E}_{\mu_m} \left[ f(Y) \right], \label{e_ace_exp} \\
&\left. \frac{d^2}{d\sigma_i d\sigma_j} A_{gce}^{f,g} \right|_{ \sigma_i , \sigma_j =0} = \frac{1}{N} \cov_{\mu^{\sigma}} \left( f(X) , g(X) \right), \label{e_agce_cov} \\
&\left. \frac{d^2}{d\sigma_i d\sigma_j} A_{ce}^{f,g} \right|_{ \sigma_i , \sigma_j =0} = \frac{1}{N} \cov_{\mu_m} \left( f(Y) , g(Y) \right), \label{e_ace_cov}
\end{align}
\end{lemma}
\medskip 

\noindent \emph{Proof of Lemma~\ref{l_second_derivative_covariance}.} \
We shall only provide the proof of~\eqref{e_agce_exp} and~\eqref{e_agce_cov}. The formula~\eqref{e_ace_exp} and~\eqref{e_ace_cov} can be derived using the same type of argument. From the definition of~$A_{gce}^{f,g}$, it holds that
\begin{align}
\frac{d}{d\sigma_i} A_{gce}^{f,g} \left(\sigma_i, \sigma_j \right) &= \frac{1}{N} \frac{ \int f(x) \exp\left( \sigma \sum_{k=1}^N x_k + \sigma_i f(x) + \sigma_j g(x) -H(x) \right) dx}{\int \exp\left( \sigma \sum_{k=1}^N x_k + \sigma_i f(x) + \sigma_j g(x) -H(x) \right) dx  } \\
& = \frac{1}{N} \mathbb{E}_{ \mu^{\sigma_i, \sigma_j } } \left[f(W) \right],
\end{align}
where~$W= \left(W_1 , W_2, \cdots, W_N \right)$ is a random variable distributed according to~$\mu^{\sigma_i, \sigma_j}$. A similar calculation yields
\begin{align}
\frac{d^2}{ d\sigma_i d\sigma_j} A_{gce}^{f,g} \left(\sigma_i, \sigma_j \right) &= \frac{1}{N}\mathbb{E}_{ \mu^{\sigma_i, \sigma_j } } \left[f(W)g(W) \right] -\frac{1}{N}\mathbb{E}_{ \mu^{\sigma_i, \sigma_j } } \left[f(W) \right]\mathbb{E}_{ \mu^{\sigma_i, \sigma_j } } \left[g(W) \right]  \\
& = \frac{1}{N} \cov_{\mu^{\sigma_i, \sigma_j } } \left( f(W), g(W) \right). 
\end{align}
Then~\eqref{e_agce_exp} and~\eqref{e_agce_cov} follow from the observation~\eqref{e_relation_between_measures}. \qed 

\medskip

Second, the Cram\'er's representation implies that the difference between two free energies is described by the density of a real-valued random variable distributed according to the modified gce~$\mu^{\sigma_i, \sigma_j}$. More precisely, we have

\begin{lemma}\label{l_Cramers_representation}
Let~$W= \left(W_1, W_2, \cdots, W_N \right)$ be a real-valued random vector distributed according to~$\mu^{\sigma_i, \sigma_j}$. Then it holds that
\begin{align}
A_{ce}^{f,g}\left(\sigma_i, \sigma_j \right) - A_{gce}^{f,g} \left(\sigma_i , \sigma_j \right) = \frac{1}{N} \ln g_{\sigma_i , \sigma_j } (0),
\end{align}
where~$g_{\sigma_i , \sigma_j }$ denotes the density of
\begin{align}
\frac{1}{\sqrt{N}} \sum_{k=1}^N \left(W_k -m\right).
\end{align}
\end{lemma}
\medskip

\noindent \emph{Proof of Lemma~\ref{l_Cramers_representation}.} \  A direct calculation yields that
\begin{align}
&A_{ce}^{f,g}\left(\sigma_i, \sigma_j \right) - A_{gce}^{f,g} \left(\sigma_i , \sigma_j \right) \\
& \qquad = \frac{1}{N} \ln \int_{\left\{\frac{1}{N} \sum_{k=1}^N x_k =m  \right\}} \exp\left( \sigma \sum_{k=1}^N x_k + \sigma_i f(x) + \sigma_j g(x) -H(x) \right) \mathcal{L}^{N-1}(dx) \\
&\qquad \quad - \frac{1}{N} \ln \int \exp\left( \sigma \sum_{k=1}^N x_k + \sigma_i f(x) + \sigma_j g(x) -H(x) \right) dx \\
& \qquad = \frac{1}{N} \ln \frac{ \int_{\left\{\frac{1}{N} \sum_{k=1}^N x_k =m  \right\}} \exp\left( \sigma \sum_{k=1}^N x_k + \sigma_i f(x) + \sigma_j g(x) -H(x) \right) \mathcal{L}^{N-1}(dx)}{\int \exp\left( \sigma \sum_{k=1}^N x_k + \sigma_i f(x) + \sigma_j g(x) -H(x) \right) dx} \\
& \qquad  =\frac{1}{N} \ln \frac{ \int_{\left\{\frac{1}{\sqrt{N}} \sum_{k=1}^N \left(x_k -m\right) =0 \right\}} \exp\left( \sigma \sum_{k=1}^N x_k + \sigma_i f(x) + \sigma_j g(x) -H(x) \right) \mathcal{L}^{N-1}(dx)}{\int \exp\left( \sigma \sum_{k=1}^N x_k + \sigma_i f(x) + \sigma_j g(x) -H(x)\right) dx} \\
& \qquad = \frac{1}{N}  \ln g_{\sigma_i , \sigma_j } (0).
\end{align}
\qed
\medskip

To derive the correct bounds for Theorem~\ref{p_equivalence_observables} and Theorem~\ref{p_equivalence_decay_of_correlations}, a careful estimation on~$g_{\sigma_i , \sigma_j}$ is needed. We begin with deducing the following representation.

\begin{lemma}\label{p_representation_local_cramer_2}
Recall the definition~\eqref{d_def_mi} of~$m_k$. It holds that
\begin{align}
2\pi  g_{0, 0 } (0) =  \int_{\mathbb{R}} \mathbb{E}_{\mu^{\sigma}} \left[ \exp\left(i \frac{1}{\sqrt{N}} \sum_{k=1}^N \left(X_k -m_k \right) \xi \right) \right] d\xi, \label{e_function_representation}
\end{align}
\begin{align}
2\pi \left.\frac{d}{d\sigma_i } g_{\sigma_i, \sigma_j } (0)\right|_{\sigma_i, \sigma_j=0} &
&= \mathbb{E}_{\mu^{\sigma}} \left[ \left(f(X) - \mathbb{E}_{\mu^{\sigma}} \left[ f(X)\right] \right) \exp\left( i \frac{1}{\sqrt{N}} \sum_{k=1}^N \left(X_k -m_k \right) \xi \right) \right], \label{e_first_derivative_representation}
\end{align}
\begin{align}
&\frac{d^2}{d\sigma_i d\sigma_j }  g_{\sigma_i, \sigma_j} (0) \bigg{|}_{\sigma_i, \sigma_j =0} \\
& \qquad = \mathbb{E}_{\mu^{\sigma}} \left[ \left(f(X) - \mathbb{E}_{\mu^{\sigma}} \left[ f(X)\right] \right)\left(g(X) - \mathbb{E}_{\mu^{\sigma}} \left[ g(X)\right] \right) \exp\left( i \frac{1}{\sqrt{N}} \sum_{k=1}^N \left(X_k -m_k \right) \xi \right) \right]. \label{e_second_derivative_representation}
\end{align}
\end{lemma}
\medskip

The proof of Lemma~\ref{p_representation_local_cramer_2} is a straight-forward application of Fourier inversion. \\

\noindent \emph{Proof of Lemma~\ref{p_representation_local_cramer_2}.} \ We start with deriving~\eqref{e_function_representation}. The inverse Fourier transform yields
\begin{align}
2\pi g_{\sigma_i, \sigma_j } (0)  &= \int_{\mathbb{R}} \mathbb{E}_{\mu^{\sigma_i, \sigma_j}} \left[ \exp\left(i \frac{1}{\sqrt{N}} \sum_{k=1}^N \left(W_k -m \right) \xi \right) \right] d\xi \\
& \overset{\eqref{e_m=sum_mi}}{=}\int_{\mathbb{R}} \mathbb{E}_{\mu^{\sigma_i, \sigma_j}} \left[ \exp\left(i \frac{1}{\sqrt{N}} \sum_{k=1}^N \left(W_k -m_k \right) \xi \right) \right] d\xi. \label{e_fourier_inversion}
\end{align}
Setting now~$\sigma_i=\sigma_j=0$ in combination with~\eqref{e_relation_between_measures} yields the desire formula~\eqref{e_function_representation}.\\

Let us now turn to the verification of~\eqref{e_first_derivative_representation}. Observe that for any smooth function~$h$,
\begin{align}
&\frac{d}{d\sigma_i} \mathbb{E}_{\mu^{\sigma_i, \sigma_j}}\left[ h(W) \right]\\
&\qquad  = \frac{d}{d\sigma_i } \int h(x) \frac{1}{Z} \exp\left( \sigma \sum_{k=1}^N x_k + \sigma_i f(x) + \sigma_j g(x) -H(x) \right)  dx \\
& \qquad = \int \frac{d}{d\sigma_i} h(x) \frac{1}{Z} \exp\left( \sigma \sum_{k=1}^N x_k + \sigma_i f(x) + \sigma_j g(x) -H(x) \right)  dx  \\
&\qquad \quad  + \int  \left(f(x) - \mathbb{E}_{\mu^{\sigma_i, \sigma_j}}\left[f(Z)\right] \right) h(x)  \frac{1}{Z}\exp\left( \sigma \sum_{k=1}^N x_k + \sigma_i x_i + \sigma_j x_j -H(x) \right)  dx \\
& \qquad = \mathbb{E}_{\mu^{\sigma_i, \sigma_j}} \left[ \frac{d}{d\sigma_i} h(W) \right] + \mathbb{E}_{\mu^{\sigma_i, \sigma_j}} \left[ \left(f(W) - \mathbb{E}_{\mu^{\sigma_i, \sigma_j}}\left[f(W)\right] \right) h(W) \right]. \label{e_derivative_formula}
\end{align}
Applying~\eqref{e_derivative_formula} to~\eqref{e_fourier_inversion} yields
\begin{align}
2\pi \frac{d}{d\sigma_i } g_{\sigma_i, \sigma_j } (0) 
& = \int_{\mathbb{R}} \frac{d}{d\sigma_i } \mathbb{E}_{\mu^{\sigma_i, \sigma_j}} \left[ \exp\left(i \frac{1}{\sqrt{N}} \sum_{k=1}^N \left(W_k -m_k \right) \xi \right) \right] d\xi \\
& = \int_{\mathbb{R}}  \mathbb{E}_{\mu^{\sigma_i, \sigma_j}} \left[\left(f(W) - \mathbb{E}_{\mu^{\sigma_i, \sigma_j}}\left[f(W)\right] \right) \exp\left(i \frac{1}{\sqrt{N}} \sum_{k=1}^N \left(W_k -m_k \right) \xi \right) \right] d\xi. \label{e_first_derivative_general}
\end{align}
Setting now~$\sigma_i=\sigma_j=0$ in combination with~\eqref{e_relation_between_measures} yields the desired formula~\eqref{e_first_derivative_representation}. The formula~\eqref{e_second_derivative_representation} also follows from the similar computations.
\qed 
\medskip

Lastly, the following Propositions (Proposition~\ref{p_equiv_observables_main_computation} and Proposition~\ref{p_exp_decay_main_computation}) provide the core estimates needed in the proof of Theorem~\ref{p_equivalence_observables} and Theorem~\ref{p_equivalence_decay_of_correlations}.

\begin{proposition} \label{p_equiv_observables_main_computation} For each~$\varepsilon>0$, there exist constants~$N_0 \in \mathbb{N}$ and~$\tilde{C}=\tilde{C}(\varepsilon) \in (0, \infty)$ independent of the external fields~$s= (s_i)$, the boundary~$x^{\mathbb{Z} \setminus [N]}$, and the mean spin~$m$ such that for all~$N \geq N_0$, it holds that
\begin{align}
&\left|\int_{\mathbb{R}} \mathbb{E}_{\mu^{\sigma}} \left[ \left(f(X) - \mathbb{E}_{\mu^{\sigma}} \left[ f(X)\right] \right) \exp\left( i \frac{1}{\sqrt{N}} \sum_{k=1}^N \left(X_k -m_k \right) \xi \right) \right] d\xi \right|\\
&\qquad \leq \tilde{C} \frac{\| \nabla f \|_{L^2 (\mu^{\sigma})} |\supp f|^{\frac{1}{2}}}{N^{ \frac{1}{2} - \varepsilon}}. \label{e_1st_derivative2} 
\end{align}
\end{proposition}
\medskip

\begin{proposition} \label{p_exp_decay_main_computation}
For each~$\varepsilon>0$, there exist constants~$N_0 \in \mathbb{N}$ and~$\tilde{C}=\tilde{C}(\varepsilon) \in (0, \infty)$ independent of the external fields~$s= (s_i)$, the boundary~$x^{\mathbb{Z} \backslash [N]}$, and the mean spin~$m$ such that for all~$N \geq N_0$, it holds that
\begin{align}
&\left| \int_{\mathbb{R}} \mathbb{E}_{\mu^{\sigma}} \left[ \left(f(X) - \mathbb{E}_{\mu^{\sigma}} \left[ f(X)\right] \right)\left(g(X) - \mathbb{E}_{\mu^{\sigma}} \left[ g(X)\right] \right) \exp\left( i \frac{1}{\sqrt{N}} \sum_{k=1}^N \left(X_k -m_k \right) \xi \right) \right] d\xi \right|\\
& \leq \tilde{C} \ C(f,g) \left( \frac{ \left( |\supp f | + |\supp g | \right)^2 }{N^{1- \varepsilon}}+ \exp\left(-C\text{dist}\left( \supp f, \supp g \right) \right)\right), \label{e_2nd_derivative}
\end{align}
where~$C(f,g)$ is given by~\eqref{e_def_c}.
\end{proposition}
\medskip

The statement of Proposition~\ref{p_equiv_observables_main_computation} and Proposition~\ref{p_exp_decay_main_computation} should be compared to Proposition~\ref{p_main_computation}. It is not surprising that the proof of Proposition~\ref{p_equiv_observables_main_computation} and Proposition~\ref{p_exp_decay_main_computation} are similar to the proof of Proposition~\ref{p_main_computation}. We prove Proposition~\ref{p_equiv_observables_main_computation} and Proposition~\ref{p_exp_decay_main_computation} in Section~\ref{s_equiv_observables_proof_of_proposition} and Section~\ref{s_exp_decay_proof_of_proposition}, respectively.\\

Let us see how these lemmas and propositions can be used to prove Theorem~\ref{p_equivalence_observables} and Theorem~\ref{p_equivalence_decay_of_correlations}. \\

\noindent \emph{Proof of Theorem~\ref{p_equivalence_observables} and Theorem~\ref{p_equivalence_decay_of_correlations}.} \ Lemma~\ref{l_second_derivative_covariance} and Lemma~\ref{l_Cramers_representation} implies that
\begin{align}
\mathbb{E}_{\mu_m} \left[ f(Y) \right] - \mathbb{E}_{\mu^{\sigma}} \left[ f(X) \right] &= \left. \frac{d}{d \sigma_i } \left( \ln g_{\sigma_i, \sigma_j} (0)\right) \right|_{ \sigma_i , \sigma_j =0}  = \frac{1}{g_{0, 0}(0)} \left(\left. \frac{d}{d \sigma_i } g_{\sigma_i, \sigma_j}(0) \right|_{ \sigma_i , \sigma_j =0}  \right)
\end{align}
Similarly, one gets
\begin{align}
&\cov_{\mu_m}\left(f(Y), g(Y) \right) - \cov_{\mu^{\sigma}}\left(f(X), g(X) \right)  \\
& \qquad= \left. \frac{d^2}{d \sigma_i d \sigma_j} \left( \ln g_{\sigma_i, \sigma_j} (0)\right) \right|_{ \sigma_i , \sigma_j =0}\\
& \qquad = - \frac{1}{\left(g_{0,0}(0) \right)^2} \left(\left.\frac{d}{d\sigma_i} g_{\sigma_i, \sigma_j}(0)\right|_{ \sigma_i , \sigma_j =0} \right)\left(\left. \frac{d}{d\sigma_j} g_{\sigma_i, \sigma_j}(0) \right|_{ \sigma_i , \sigma_j =0} \right)\\
& \qquad \quad + \frac{1}{g_{0,0}(0)}\left(\left. \frac{d^2}{d \sigma_i d \sigma_j} g_{\sigma_i, \sigma_j}(0) \right|_{ \sigma_i , \sigma_j =0} \right).
\end{align}
Then a combination of Lemma~\ref{p_representation_local_cramer_2}, Proposition~\ref{p_main_computation}, Proposition~\ref{p_equiv_observables_main_computation} and Proposition~\ref{p_exp_decay_main_computation} provides the desired bounds of Theorem~\ref{p_equivalence_observables} and Theorem~\ref{p_equivalence_decay_of_correlations}. \qed

\medskip

Let us now give the proof of Theorem~\ref{p_exponential_decay_of_correlations_of_canonical_ensemble}. \\

\noindent \emph{Proof of Theorem~\ref{p_exponential_decay_of_correlations_of_canonical_ensemble}.} \ The desired statement is a direct consequence of Theorem~\ref{p_decay_of_correlations_gce} and Theorem~\ref{p_equivalence_decay_of_correlations}. Indeed, it holds that
\begin{align}
\left| \cov_{\mu_m ^{\Lambda}} \left( f, g \right)\right| &\leq \left|  \cov_{\mu_m^{\Lambda}}  (f,g) - \cov_{\mu^{\Lambda, \sigma}}  (f,g)  \right| + \left| \cov_{\mu^{\Lambda, \sigma}}  (f,g) \right| \\
& \leq \tilde{C} \ C(f,g) \left( \frac{ \left( |\supp f | + |\supp g | \right)^2 }{|\Lambda|^{1- \varepsilon}} + \exp\left(-C\text{dist}\left( \supp f, \supp g \right) \right) \right)  \\
 & \quad + C \|\nabla f\|_{L^2 (\mu^{\Lambda, \sigma})}\|\nabla g\|_{L^2 (\mu^{\Lambda, \sigma})} \exp \left( -C \text{dist} \left( \supp f , \supp g \right) \right) \\
 & \lesssim \tilde{C} \ C(f,g) \left( \frac{ \left( |\supp f | + |\supp g | \right)^2 }{|\Lambda|^{1- \varepsilon}} + \exp\left(-C\text{dist}\left( \supp f, \supp g \right) \right) \right) .
\end{align}
\qed


\subsection{Proof of Proposition~\ref{p_equiv_observables_main_computation}.} \label{s_equiv_observables_proof_of_proposition}
We may assume without loss of generality that
\begin{align}
\mathbb{E}_{\mu^{\sigma}} \left[ f(X) \right] = 0.
\end{align}
Let us begin with dividing the integral into inner and outer parts. We fix~$\delta>0$ small enough and decompose the integral as follows:
\begin{align}
&\int_{\mathbb{R}} \mathbb{E}_{\mu ^{\sigma} } \left[ f(X) \exp\left( i\frac{1}{\sqrt{N}} \sum_{k=1}^N \left(X_k -m_k \right) \xi \right) \right] d\xi \\
& \qquad =\int_{\{ \left| \left(1/ \sqrt{N}\right)\xi\right|\leq\delta\}} \mathbb{E}_{\mu^{\sigma}  } \left[ f(X) \exp\left( i\frac{1}{\sqrt{N}} \sum_{k=1}^N \left(X_k -m_k \right) \xi \right) \right] d\xi \label{e_equiv_observables_inner}  \\
& \qquad \quad + \int_{\{ \left| \left(1/ \sqrt{N}\right)\xi\right| > \delta\}} \mathbb{E}_{\mu^{\sigma} } \left[ f(X) \exp\left( i\frac{1}{\sqrt{N}} \sum_{k=1}^N \left(X_k -m_k \right) \xi \right) \right] d\xi. \label{e_equiv_observables_outer}
\end{align}
Our aim is to bound the inner and outer integrals separately. More precisely, we prove the following two statements:

\begin{lemma} [Estimation of the inner integral] \label{l_equiv_observables_inner}
Fix~$\delta>0$ small enough. For each~$\varepsilon >0$, there exist constants~$N_1  \in \mathbb{N}$ and~$\tilde{C}= \tilde{C} (\epsilon) \in (0, \infty) $ independent of the external fields~$s=(s_i)$, the boundary~$x^{\mathbb{Z} \setminus [N]}$, and the mean spin~$m$ such that for any~$N \geq N_1$,
\begin{align}
&\left| \int_{\{ \left| \left(1/ \sqrt{N}\right)\xi\right|\leq\delta\}} \mathbb{E}_{\mu^{\sigma}  } \left[ f(X) \exp\left( i\frac{1}{\sqrt{N}} \sum_{k=1}^N \left(X_k -m_k \right) \xi \right) \right] d\xi \right|  \leq \tilde{C} \frac{ \| \nabla f \|_{L^2 (\mu^{\sigma})} |\supp f |^{\frac{1}{2}} } {N^{ \frac{1}{2} - \varepsilon}}.
\end{align}
\end{lemma}
\medskip 

\begin{lemma} [Estimation of the outer integral] \label{l_equiv_observables_outer}
Fix~$\delta>0$ small enough. For each~$\varepsilon >0$, there exist constants~$N_2  \in \mathbb{N}$ and~$\tilde{C}= \tilde{C} (\epsilon) \in (0, \infty) $ independent of the external fields~$s=(s_i)$, the boundary~$x^{\mathbb{Z} \setminus [N]}$, and the mean spin~$m$ such that for any~$N \geq N_2$,
\begin{align}
&\left| \int_{\{ \left| \left(1/ \sqrt{N}\right)\xi\right|>\delta\}} \mathbb{E}_{\mu^{\sigma}  } \left[ f(X) \exp\left( i\frac{1}{\sqrt{N}} \sum_{k=1}^N \left(X_k -m_k \right) \xi \right) \right] d\xi \right|  \leq \tilde{C} \frac{ \| \nabla f \|_{L^2 (\mu^{\sigma})} |\supp f |^{\frac{1}{2}} } {N^{ \frac{1}{2} - \varepsilon}}.
\end{align}
\end{lemma}
\medskip

The proof of Proposition~\ref{p_equiv_observables_main_computation} is a direct consequence of Lemma~\ref{l_equiv_observables_inner} and Lemma~\ref{l_equiv_observables_outer}. \\

\noindent \emph{Proof of Proposition~\ref{p_equiv_observables_main_computation}.} \ For fixed~$\delta>0$ small enough, we choose~$N_1, N_2 \in \mathbb{N}$ as in Lemma~\ref{l_equiv_observables_inner} and Lemma~\ref{l_equiv_observables_outer}. Define~$N_0 = \max (N_1, N_2)$. Then for any~$N \geq N_0$, the triangle inequality followed by Lemma~\ref{l_equiv_observables_inner} and Lemma~\ref{l_equiv_observables_outer} yields
\begin{align}
&\left| \int_{\mathbb{R}} \mathbb{E}_{\mu ^{\sigma} } \left[ f(X) \exp\left( i\frac{1}{\sqrt{N}} \sum_{k=1}^N \left(X_k -m_k \right) \xi \right) \right] d\xi \right| \\
& \qquad \leq \left| \int_{\{ \left| \left(1/ \sqrt{N}\right)\xi\right|\leq\delta\}} \mathbb{E}_{\mu^{\sigma}  } \left[ f(X) \exp\left( i\frac{1}{\sqrt{N}} \sum_{k=1}^N \left(X_k -m_k \right) \xi \right) \right] d\xi \right| \\
& \qquad \quad + \left| \int_{\{ \left| \left(1/ \sqrt{N}\right)\xi\right|>\delta\}} \mathbb{E}_{\mu^{\sigma}  } \left[ f(X) \exp\left( i\frac{1}{\sqrt{N}} \sum_{k=1}^N \left(X_k -m_k \right) \xi \right) \right] d\xi \right| \\
& \qquad \lesssim \frac{ \| \nabla f \|_{L^2 (\mu^{\sigma})} |\supp f |^{\frac{1}{2}} } {N^{ \frac{1}{2} - \varepsilon}}.
\end{align}
\qed
\medskip

The rest of the section consists of the proof of Lemma~\ref{l_equiv_observables_inner} and Lemma~\ref{l_equiv_observables_outer}, which contain lengthy calculations. \\

\begin{figure}[t]
\centering
\begin{tikzpicture}[xscale=1.3]

\draw[fill] (0,0) circle [radius=0.05];
\draw[fill] (.25,0) circle [radius=0.05];
\draw[fill] (.5,0) circle [radius=0.05];
\draw[fill] (.75,0) circle [radius=0.05];
\draw[fill] (1,0) circle [radius=0.05];

\draw[fill] (1.475,0) circle [radius=0.02];
\draw[fill] (1.625,0) circle [radius=0.02];
\draw[fill] (1.775,0) circle [radius=0.02];

\draw[fill] (2.25,0) circle [radius=0.05];
\draw[fill] (2.5,0) circle [radius=0.05];
\draw[fill] (2.75,0) circle [radius=0.05];
\draw[fill] (3,0) circle [radius=0.05];
\draw[fill] (3.25,0) circle [radius=0.05];

\draw (3.5,0) circle[radius=0.05];
\draw (3.75,0) circle[radius=0.05];
\draw (4,0) circle[radius=0.05];
\draw (4.25,0) circle[radius=0.05];
\draw (4.5,0) circle[radius=0.05];
\draw (4.75,0) circle[radius=0.05];
\draw (5,0) circle[radius=0.05];
\draw (5.25,0) circle[radius=0.05];
\draw (5.5,0) circle[radius=0.05];

\draw[fill] (5.75,0) circle [radius=0.05];
\draw[fill] (6,0) circle [radius=0.05];
\draw[fill] (6.25,0) circle [radius=0.05];
\draw[fill] (6.5,0) circle [radius=0.05];
\draw[fill] (6.75,0) circle [radius=0.05];

\draw[fill] (7.225,0) circle [radius=0.02];
\draw[fill] (7.375,0) circle [radius=0.02];
\draw[fill] (7.525,0) circle [radius=0.02];

\draw[fill] (8,0) circle [radius=0.05];
\draw[fill] (8.25,0) circle [radius=0.05];
\draw[fill] (8.5,0) circle [radius=0.05];
\draw[fill] (8.75,0) circle [radius=0.05];
\draw[fill] (9,0) circle [radius=0.05];
\node[align=center, below] at (0,-.05) {$1$};
\node[align=center, below] at (4.5,-.05) {$l$};
\node[align=center, below] at (9,-.05) {$N$};
\node[align=center, below] at (4.5,-1.2) {};

\draw[decorate,decoration={brace,mirror},thick] (0,-.5) -- node[below]{$E_{2}^{f}$} (3.25,-.5);
\draw[decorate,decoration={brace,mirror},thick] (3.5,-.5) -- node[below]{$E_{1}^{f}$} (5.5,-.5);
\draw[decorate,decoration={brace,mirror},thick] (5.75,-.5) -- node[below]{$E_{2}^{f}$} (9,-.5);

\draw[thick,<->] (3.5,.5) -- (4.45,.5);
\node[align=center, above] at (4,.5) {$L$};
\draw[thick,<->] (4.55,.5) -- (5.5,.5);
\node[align=center, above] at (5,.5) {$L$};
\end{tikzpicture}
\caption{The sets $E_1^{f}$ and $E_2 ^{f}$ where~$f(x)=x_i-m_i$}\label{f_sets_E}
\end{figure}
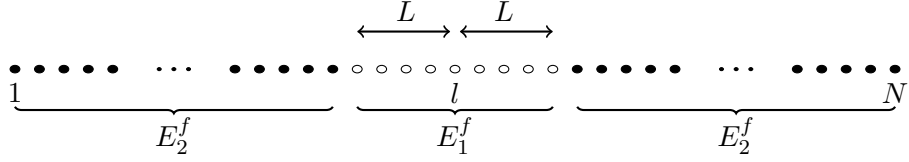

Let us begin with providing auxiliary ingredients for deducing Lemma~\ref{l_equiv_observables_inner}. We define the auxiliary sets~$E_1^{f}$ and~$E_2^{f}$ as (see Figure~\ref{f_sets_E}):
\begin{align}
&E_1^{f} : = \{1, 2, \cdots, N \} \cap \{ k : \ \text{dist}\left( k, \supp(f)\right) \leq L  \}, \\
&E_2^{f} : = \{1, 2, \cdots, N \} \cap \{ k : \ \text{dist}\left( k, \supp(f)\right) > L  \},
\end{align}
where~$L\ll N$ is a positive integer that will be chosen later. The main ingredients for the proof of Lemma~\ref{l_equiv_observables_inner} are Theorem~\ref{p_decay_of_correlations_gce}, Lemma~\ref{l_amgm} and an extension of~\cite[Lemma 7]{KwMe18}.

\begin{lemma}[Extension of Lemma 7 in~\cite{KwMe18}] \label{l_observables_exchanging_exponential_terms}
For large enough~$N$ and~$\delta>0$ small enough, there exists a positive constant~$C>0$ such that the following inequalities hold for all~$\xi \in \mathbb{R}$ with~$ \frac{\left|\xi\right|}{\sqrt{N}}  \leq \delta$.
\begin{align}
 \left|\mathbb{E}_{\mu^{\sigma}}\left[\exp\left( i \sum_{k \in E_2^{f} } \left( X_k -m_k \right) \frac{\xi}{\sqrt{N}}\right) \mid \mathscr{F}_{f} \right] \right| \leq C \left( 1 + \xi^2 \right)\exp\left(-C\xi^2 \right) .
\label{e_lemma exchanging 2nd der}
\end{align}
where~$\mathscr{F}_{f}$ denotes the sigma algebra defined by
\begin{align}
\mathscr{F}_{f} : = \sigma\left(X_k , \ k \in E_1 ^{f} \right).
\end{align} 
\end{lemma}
\medskip 

\begin{remark} \label{r_exchanging_exponential_terms}
The proof of Lemma~\ref{l_observables_exchanging_exponential_terms} is almost similar to that of~\cite[Lemma 7]{KwMe18}. One should compare the sets~$E_1 ^{f}, E_2^{f}$ with~$E_1^{l}$ and~$E_2^{l}$ in~\cite{KwMe18}. The main difference is that we assume finite range interaction with range~$R$ instead of the nearest neighbor interaction. However, there is only a cosmetic difference between these two proofs. We leave the details to the reader. 
\end{remark}
\medskip

\noindent \emph{Proof of Lemma~\ref{l_equiv_observables_inner}.} \ Let us denote~$e : \mathbb{R}^2 \to \mathbb{C}$ by
\begin{align}
e(\xi_1, \xi_2 ) = \mathbb{E}_{\mu^{\sigma}} \left[f(X) \exp \left( i  \sum_{k \in E_1 ^{f}}\left( X_k - m_k \right) \xi_1 + i  \sum_{l \in E_2 ^{f}} \left( X_l -m_l \right) \xi_2 \right) \right].
\end{align}
Then a Taylor expansion with respect to the first variable~$\xi_1$ yields
\begin{align}
e(\xi_1 , \xi_2 ) &=  \mathbb{E}_{\mu^{\sigma}} \left[f(X) \exp \left(  i  \sum_{l \in E_2 ^{f}} \left( X_l -m_l \right) \xi_2 \right) \right]   \\
& \quad + \mathbb{E}_{\mu^{\sigma}} \left[f(X)\left( i \sum_{k_1 \in E_1 ^{f}} \left( X_{k_1} - m_{k_1} \right) \right) \right. \\
& \left. \qquad \qquad \qquad \times\exp \left( i  \sum_{k_2 \in E_1 ^{f}}\left( X_{k_2} - m_{k_2} \right) \tilde{\xi}_1 + i  \sum_{l \in E_2 ^{f}} \left( X_l -m_l \right) \xi_2 \right) \right]\xi_1 ,
\end{align}
where~$\tilde{\xi}_1$ is a real number between~$0$ and~$\xi_1$. In particular for~$(\xi_1, \xi_2 ) = \left( \frac{\xi}{\sqrt{N}}, \frac{\xi}{\sqrt{N}} \right)$ it holds that
\begin{align}
&\mathbb{E}_{\mu^{\sigma}} \left[f(X) \exp \left( i \frac{1}{\sqrt{N}} \sum_{k \in E_1 ^{f}}\left( X_k - m_k \right) \xi + i \frac{1}{\sqrt{N}} \sum_{l \in E_2 ^{f}} \left( X_l -m_l \right) \xi \right) \right] \label{e_equiv_observables_integrand}\\
& \qquad = \mathbb{E}_{\mu^{\sigma}} \left[f(X) \exp \left(  i \frac{1}{\sqrt{N}} \sum_{l \in E_2 ^{f}} \left( X_l -m_l \right) \xi \right) \right]  \label{e_equiv_observables_0th_taylor}  \\
& \qquad \quad + \mathbb{E}_{\mu^{\sigma}} \left[f(X)\left(   \sum_{k_1 \in E_1 ^{f}} \left( X_{k_1} - m_{k_1} \right) \right) \right. \\
& \left. \qquad \qquad \qquad \quad \times\exp \left( i \frac{1}{\sqrt{N}} \sum_{k_2 \in E_1 ^{f}}\left( X_{k_2} - m_{k_2} \right) \frac{\tilde{\xi}}{\sqrt{N}} + i  \sum_{l \in E_2 ^{f}} \left( X_l -m_l \right) \frac{\xi}{\sqrt{N}} \right) \right] i \frac{\xi}{\sqrt{N}}, \label{e_equiv_observables_1st_taylor}
\end{align}
where~$ \frac{\tilde{\xi}}{\sqrt{N}}$ is a real number between~$0$ and~$\frac{\xi}{\sqrt{N}}$. Then it holds from Theorem~\ref{p_decay_of_correlations_gce} that
\begin{align}
\left|T_{\eqref{e_equiv_observables_0th_taylor}}\right| & = \left| \cov_{\mu^{\sigma}} \left( f(X) , \exp \left( i \frac{1}{\sqrt{N}} \sum_{l \in E_2 ^f} \left( X_l -m_l \right) \xi \right) \right) \right| \\
&\lesssim \| \nabla f \|_{L^2 (\mu^{\sigma})} |\xi| \exp \left( - CL \right).
\end{align}
Let us turn to the estimation of~\eqref{e_equiv_observables_1st_taylor}. Considering the conditional expectation with respect to the sigma algebra~$\mathscr{F}_{f} = \sigma \left( X_k , \  k \in E_1^{f} \right)$, we write
\begin{align}
T_{\eqref{e_equiv_observables_1st_taylor}}   &=  \mathbb{E}_{\mu^{\sigma}} \left[ f(X) \left(  \sum_{k_1 \in E_1 ^{f} } \left(X_{k_1} -m_{k_1} \right) \right)   \exp\left( i \sum_{k_2 \in E_1 ^{f}} \left(X_{k_2} -m_{k_2} \right) \frac{\xi}{\sqrt{N}} \right) \right.  \\
&  \left. \qquad \qquad \qquad  \times \mathbb{E}_{\mu^{\sigma}} \left[ \exp\left(  \sum_{l \in E_2 ^{f}} \left(X_l -m_l \right) \frac{\xi}{\sqrt{N}} \right) \mid \mathscr{F}_{f,g}  \right] \right]  i\frac{\xi}{\sqrt{N}}.
\end{align}
Now we apply Lemma~\ref{l_observables_exchanging_exponential_terms} followed by Cauchy-Schwartz inequality, Poincar\'e inequality and Lemma~\ref{l_variance_estimate} to obtain
\begin{align}
\left|T_{\eqref{e_equiv_observables_1st_taylor}}\right| &\leq \mathbb{E}_{\mu^{\sigma}} \left| f(X) \left( \sum_{k_1 \in E_1 ^{f}} \left( X_{k_1} - m_{k_1} \right) \right) \right| (1+\xi^2) \exp\left( - C \xi^2 \right) \frac{|\xi|}{\sqrt{N}} \\
& \leq \left(\mathbb{E}_{\mu^{\sigma}} \left[ f(X) ^2 \right] \right)^{\frac{1}{2}} \left( \mathbb{E}_{\mu^{\sigma}} \left[ \left(\sum_{k_1 \in E_1 ^{f}} \left( X_{k_1} - m_{k_1} \right) \right)^2 \right] \right)^{\frac{1}{2}} (1+\xi^2 ) \exp \left( - C\xi^2 \right) \frac{|\xi|}{\sqrt{N}} \\
& \lesssim \| \nabla f \|_{L^2 (\mu^{\sigma})} \left(\var_{\mu^{\sigma}} \left( \sum_{k_1 \in E_1 ^{f}} X_{k_1} \right) \right)^{\frac{1}{2}} (1+\xi^2 ) \exp \left( - C\xi^2 \right) \frac{|\xi|}{\sqrt{N}} \\
& \lesssim \| \nabla f \|_{L^2 (\mu^{\sigma})} | E_1 ^{f}| ^{\frac{1}{2}} (1+\xi^2 ) \exp \left( - C\xi^2 \right) \frac{|\xi|}{\sqrt{N}}.
\end{align}
Note that cardinality~$|E_1 ^{f}|$ is bounded from above by~$|\supp f | \cdot 2L$. Therefore we conclude that
\begin{align}
\left|T_{\eqref{e_equiv_observables_integrand}}\right| & \leq \left|T_{\eqref{e_equiv_observables_0th_taylor}} \right| +\left|T_{\eqref{e_equiv_observables_1st_taylor}} \right| \\
& \lesssim \| \nabla f \|_{L^2 (\mu^{\sigma})} |\xi| \exp \left( - CL \right) + \| \nabla f \|_{L^2 (\mu^{\sigma})} | \supp f| ^{\frac{1}{2}} (1+\xi^2 ) \exp \left( - C\xi^2 \right)|\xi| \frac{L^{\frac{1}{2}}}{\sqrt{N}}.
\end{align}
Note that for~$L=N^{\varepsilon} \ll N$ and~$N$ large enough, it holds that
\begin{align}
\int_{\{ \left| \left(1/ \sqrt{N}\right)\xi\right|\leq\delta\}} |\xi| \exp\left(-CL\right) d\xi &\lesssim \frac{1}{\sqrt{N}}, \\
\int_{\{ \left| \left(1/ \sqrt{N}\right)\xi\right|\leq\delta\}} |\xi|(1+ \xi^2) \exp\left(-C\xi^2\right) d\xi &\lesssim 1 
\end{align}
Hence we have the desired bound
\begin{align}
&\left| \int_{\{ \left| \left(1/ \sqrt{N}\right)\xi\right|\leq\delta\}} \mathbb{E}_{\mu^{\sigma}} \left[f(X) \exp \left( i \frac{1}{\sqrt{N}} \sum_{k \in E_1 ^{f}}\left( X_k - m_k \right) \xi + i \frac{1}{\sqrt{N}} \sum_{l \in E_2 ^{f}} \left( X_l -m_l \right) \xi \right) \right] d\xi \right| \\
& \qquad \lesssim \| \nabla f \|_{L^2 (\mu^{\sigma})} \frac{1}{\sqrt{N}} + \| \nabla f \|_{L^2 (\mu^{\sigma})} |\supp f |^{\frac{1}{2}} \frac{1}{N^{\frac{1}{2} - \varepsilon}}\\
& \qquad \lesssim \| \nabla f \|_{L^2 (\mu^{\sigma})} |\supp f |^{\frac{1}{2}} \frac{1}{N^{\frac{1}{2} - \varepsilon}}.
\end{align}
\qed
\medskip

Let us turn to the proof of Lemma~\ref{l_equiv_observables_outer}. This part is motivated by the argument presented in~\cite{KwMe18}. The main difference is again, we consider the finite range interaction with range~$R$ instead of the nearest-neighbor interaction. \\

Consider the characteristic function~$\varphi_W (\xi)$ of the random variable
\begin{align}
W= \frac{1}{\sqrt{N}}\sum_{k=1}^N \left( X_k -m_k \right)    
\end{align}
which is given by
\begin{align}
\varphi_W (\xi) = \mathbb{E}_{\mu^{\sigma}} \left[ \exp \left( i \frac{1}{\sqrt{N}}\sum_{k=1}^N \left( X_k -m_k \right)\xi \right)\right].
\end{align}

\begin{figure}[t]
\centering
\begin{tikzpicture}[xscale=1.3]

\draw (.5,0) node[cross] {};
\draw (.75,0) node[cross] {};
\draw (1,0) node[cross] {};
\draw (1.25,0) node[cross] {};
\draw (1.5,0) node[cross] {};
\draw (1.75,0) node[cross] {};
\draw[fill] (2,0) circle[radius=0.05];

\draw (2.25,0) node[cross] {};
\draw (2.5,0) node[cross] {};
\draw (2.75,0) node[cross] {};
\draw (3,0) node[cross] {};
\draw (3.25,0) node[cross] {};
\draw (3.5,0) node[cross] {};
\draw[fill] (3.75,0) circle[radius=0.05];

\draw (4,0) node[cross] {};
\draw (4.25,0) node[cross] {};
\draw (4.5,0) node[cross] {};
\draw (4.75,0) node[cross] {};
\draw (5,0) node[cross] {};
\draw (5.25,0) node[cross] {};
\draw[fill] (5.5,0) circle[radius=0.05];

\draw (5.75,0) node[cross] {};
\draw (6,0) node[cross] {};
\draw (6.25,0) node[cross] {};
\draw (6.5,0) node[cross] {};
\draw (6.75,0) node[cross] {};
\draw (7,0) node[cross] {};
\draw[fill] (7.25,0) circle[radius=0.05];

\draw (7.5,0) node[cross] {};
\draw (7.75,0) node[cross] {};
\draw (8,0) node[cross] {};
\draw (8.25,0) node[cross] {};
\draw (8.5,0) node[cross] {};
\draw (8.75,0) node[cross] {};

\draw[fill] (9.05,0) circle [radius=0.02];
\draw[fill] (9.15,0) circle [radius=0.02];
\draw[fill] (9.25,0) circle [radius=0.02];

\node[align=center, below] at (2,-.05) {$R+1$};
\node[align=center, below] at (3.75,-.05) {$2(R+1)$};
\node[align=center, below] at (5.5,-.05) {$3(R+1)$};
\node[align=center, below] at (7.25,-.05) {$4(R+1)$};

\draw[thick,<->] (.5,.5) -- (1.75,.5);
\node[align=center, above] at (1.125,.5) {$R$};
\draw[thick,<->] (2.25,.5) -- (3.5,.5);
\node[align=center, above] at (2.875,.5) {$R$};
\draw[thick,<->] (4,.5) -- (5.25,.5);
\node[align=center, above] at (4.625,.5) {$R$};
\draw[thick,<->] (5.75,.5) -- (7,.5);
\node[align=center, above] at (6.375,.5) {$R$};
\draw[thick,<->] (7.5,.5) -- (8.75,.5);
\node[align=center, above] at (8.125,.5) {$R$};

\end{tikzpicture}
\caption{The set~$[N]_R$ where~$R=6$.  }\label{f_conditioned_rv}
\end{figure}
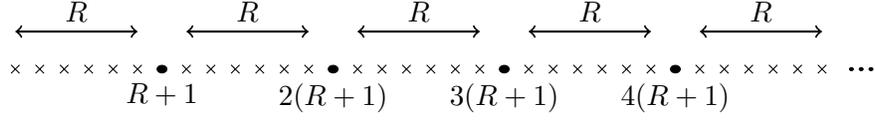

Our strategy is to induce an artificial independence by conditioning. More precisely, let $[N]_R : = \left\{R+1, 2(R+1), \cdots \right\} \cap [N] $ (see Figure~\ref{f_conditioned_rv}).
Because of the finite range interaction with range~$R$ we have a product structure of conditional characteristic functions i.e. 

\newpage
\begin{align}
&\mathbb{E}_{\mu^{\sigma}} \left[ \exp \left( i \frac{1}{\sqrt{N}}\sum_{k=1}^N \left( X_k -m_k \right)\xi \right)\right] \\
& \qquad =\mathbb{E}_{\mu^{\sigma}} \left[\exp \left( i \frac{1}{\sqrt{N}}\sum_{k \in [N] \backslash [N]_R} \left( X_k -m_k \right)\xi \right) \right. \\
& \qquad \qquad \qquad  \left. \times \mathbb{E}_{\mu^{\sigma}} \left[ \exp \left( i \frac{1}{\sqrt{N}}\sum_{l \in  [N]_R} \left( X_l -m_l \right)\xi \right) \Big| \  X_j, j \in [N] \backslash [N]_R \right] \right] \\
& \qquad = \mathbb{E}_{\mu^{\sigma}} \left[\exp \left( i \frac{1}{\sqrt{N}}\sum_{k \in [N] \backslash [N]_R} \left( X_k -m_k \right)\xi \right) \right. \\
& \qquad \qquad \qquad  \left. \times \prod_{l \in [N]_R} \mathbb{E}_{\mu^{\sigma}} \left[ \exp \left( i \frac{1}{\sqrt{N}} \left( X_l -m_l \right)\xi \right)  \Big| \  X_j, j \in [N] \backslash [N]_R  \right] \right].
\end{align}
The key ingredient for deducing Lemma~\ref{l_equiv_observables_outer} is the following statement:

\begin{lemma}[Extension of (76) and (77) in~\cite{KwMe18}] \label{l_outer_key_ingredient} 
For any~$\delta>0$ small, there is $\lambda =\lambda(\delta) \in (0, 1)$ such that for any~$l \in [N]_R$
\begin{align}
\left|\mathbb{E}_{\mu^{\sigma}} \left[ \exp \left( i \frac{1}{\sqrt{N}} \left( X_l -m_l \right)\xi \right)  \Big| \  X_j, j \in [N] \backslash [N]_R  \right]\right| \leq \lambda \qquad \text{for all } |\xi| > \delta \label{e_outer_key_ingredient_1}
\end{align}
and
\begin{align}
\left|\mathbb{E}_{\mu^{\sigma}} \left[ \exp \left( i \frac{1}{\sqrt{N}} \left( X_l -m_l \right)\xi \right)  \Big| \  X_j, j \in [N] \backslash [N]_R  \right]\right| \lesssim \frac{1}{|\xi|}.\label{e_outer_key_ingredient_2}
\end{align} 
\end{lemma}
\medskip

We refer to~\cite{KwMe18} or~\cite{GrOtViWe09} for the proof of Lemma~\ref{l_outer_key_ingredient}. Now we are ready to give the proof of Lemma~\ref{l_equiv_observables_outer}. \\

\noindent \emph{Proof of Lemma~\ref{l_equiv_observables_outer}.} \ For each~$k \in E_2 ^f$ denote
\begin{align}
\bar{m}_{k} : = \mathbb{E}_{\mu^{\sigma}} \left[ X_k \Big| \ X_j, j \in E_1 ^f \right].
\end{align}
We write the integrand in the statement as
\begin{align}
&\mathbb{E}_{\mu^{\sigma}} \left[ f(X) \exp \left( i \frac{1}{\sqrt{N}} \sum_{k=1}^{N} \left( X_k -m_k \right) \xi  \right) \right] \\
& \qquad = \mathbb{E}_{\mu^{\sigma}} \left[  f(X) \exp \left( i \frac{1}{\sqrt{N}} \sum_{k \in E_1 ^f} \left( X_k -m_k \right) \xi + i \frac{1}{\sqrt{N}}\sum_{k \in E_2 ^f} \left(\bar{m}_k - m_k \right)\xi   \right) \right. \\
&\left. \qquad \qquad \qquad \qquad \qquad \qquad \times \mathbb{E}_{\mu^{\sigma}} \left[ \exp \left(  i \frac{1}{\sqrt{N}} \sum_{k \in E_2 ^f} \left( X_k -\bar{m}_k \right) \xi \right) \Big| \ X_j, j \in E_1 ^f \right] \right].
\end{align}
Adapting the idea presented before Lemma~\ref{l_outer_key_ingredient}, we write
\begin{align}
&\mathbb{E}_{\mu^{\sigma}} \left[ \exp \left(  i \frac{1}{\sqrt{N}} \sum_{k \in E_2 ^f} \left( X_k -\bar{m}_k \right) \xi \right) \Big| \ X_j, j \in E_1 ^f \right] \\
&  = \mathbb{E}_{\mu^{\sigma}} \left[\exp \left( i \frac{1}{\sqrt{N}}\sum_{k \in E_2 ^f \backslash [N]_R} \left( X_k -\bar{m}_k \right)\xi \right) \right. \\
& \qquad \qquad  \left. \times \prod_{l \in E_2 ^f \cap [N]_R} \mathbb{E}_{\mu^{\sigma}} \left[ \exp \left( i \frac{1}{\sqrt{N}} \left( X_l -\bar{m}_l \right)\xi \right)  \Big| \  X_k, k \in [N] \backslash [N]_R  \right] \Big| \ X_j, j \in E_1 ^f \right].
\end{align}
Denote the cardinality of the set~$E_2 ^f \cap [N]_R$ by~$N_R ^f$. We apply~\eqref{e_outer_key_ingredient_1} (on~$N_R ^f -2$ factors of~$N_R ^f$ factors) and~\eqref{e_outer_key_ingredient_2} (on the remaining~$2$ factors) to obtain
\begin{align}
\left|\mathbb{E}_{\mu^{\sigma}} \left[ \exp \left(  i \frac{1}{\sqrt{N}} \sum_{k \notin S} \left( X_k -\bar{m}_k \right) \xi \right) \Big| \ X_j, j \in S \right]\right| &\lesssim  \lambda^{N_R ^f-2} \left( \frac{1}{1 + \left( 1/ \sqrt{N}\right)|\xi| }\right)^2 \\
& \lesssim N\lambda^{N_R ^f-2} \frac{1}{1+ \xi^2},
\end{align}
and as a consequence,
\begin{align}
\left| \mathbb{E}_{\mu^{\sigma}} \left[ f(X) \exp \left( i \frac{1}{\sqrt{N}} \sum_{k=1}^{N} \left( X_k -m_k \right) \xi  \right) \right] \right| 
&  \lesssim \mathbb{E}_{\mu^{\sigma}}\left[ |f(X)| \right]N\lambda^{N_R ^f-2} \frac{1}{1+ \xi^2} \\
& \overset{Poincare}{\lesssim} \| \nabla f \|_{L^2(\mu^{\sigma})} N\lambda^{N_R ^f -2} \frac{1}{1+ \xi^2}.
\end{align}
Now we conclude that
\begin{align}
&\left| \int_{\{ \left| \left(1/ \sqrt{N}\right)\xi\right| > \delta\}} \mathbb{E}_{\mu^{\sigma} } \left[f(X) \exp\left( i\frac{1}{\sqrt{N}} \sum_{k=1}^N \left(X_k -m_k \right) \xi \right) \right] d\xi\right| \\
& \qquad \lesssim \int_{\{ \left| \left(1/ \sqrt{N}\right)\xi\right| > \delta\}} \| \nabla f \|_{L^2(\mu^{\sigma})} N\lambda^{N_R ^f-2} \frac{1}{1+ \xi^2} d\xi \\
& \qquad \lesssim \| \nabla f \|_{L^2(\mu^{\sigma})} N\lambda^{N_R ^f -2}  \overset{\lambda < 1 \ll N}{\lesssim} \frac{ \| \nabla f \|_{L^2 (\mu^{\sigma})} |\supp f |^{\frac{1}{2}} } {N^{ \frac{1}{2} - \varepsilon}}.
\end{align}
\qed

\bigskip

\subsection{Proof of Proposition~\ref{p_exp_decay_main_computation}} \label{s_exp_decay_proof_of_proposition}

The proof of Proposition~\ref{p_exp_decay_main_computation} is similar to that of Proposition~\ref{p_equiv_observables_main_computation}. The main difference is that we apply the second order Taylor expansion when estimating the inner integral. \\

As in Section~\ref{s_equiv_observables_proof_of_proposition}, we assume without loss of generality that
\begin{align}
\mathbb{E}_{\mu^{\sigma}} \left[ f(X) \right] = \mathbb{E}_{\mu^{\sigma}} \left[ g(X) \right]  =0.
\end{align}
We then divide the integral into inner and outer parts and estimate them separately. Let us fix~$\delta>0$ small enough and decompose the integral as follows:
\begin{align}
&\int_{\mathbb{R}} \mathbb{E}_{\mu ^{\sigma} } \left[ f(X) g(X) \exp\left( i\frac{1}{\sqrt{N}} \sum_{k=1}^N \left(X_k -m_k \right) \xi \right) \right] d\xi \\
& \qquad =\int_{\{ \left| \left(1/ \sqrt{N}\right)\xi\right|\leq\delta\}} \mathbb{E}_{\mu^{\sigma}  } \left[ f(X) g(X) \exp\left( i\frac{1}{\sqrt{N}} \sum_{k=1}^N \left(X_k -m_k \right) \xi \right) \right] d\xi \label{e_exp_decay_inner}  \\
& \qquad \quad + \int_{\{ \left| \left(1/ \sqrt{N}\right)\xi\right| > \delta\}} \mathbb{E}_{\mu^{\sigma} } \left[ f(X) g(X) \exp\left( i\frac{1}{\sqrt{N}} \sum_{k=1}^N \left(X_k -m_k \right) \xi \right) \right] d\xi. \label{e_exp_decay_outer}
\end{align}

The following two statements are analogues of Lemma~\ref{l_equiv_observables_inner} and Lemma~\ref{l_equiv_observables_outer}.

\begin{lemma} [Estimation of the inner integral] \label{l_exp_decay_inner}
Fix~$\delta>0$ small enough. For each~$\varepsilon >0$, there exist constants~$N_1  \in \mathbb{N}$ and~$\tilde{C}= \tilde{C} (\epsilon) \in (0, \infty) $ independent of the external fields~$s=(s_i)$, the boundary~$x^{\mathbb{Z} \setminus [N]}$, and the mean spin~$m$ such that for any~$N \geq N_1$,
\begin{align}
&\left| \int_{\{ \left| \left(1/ \sqrt{N}\right)\xi\right|\leq\delta\}} \mathbb{E}_{\mu^{\sigma}  } \left[ f(X)g(X) \exp\left( i\frac{1}{\sqrt{N}} \sum_{k=1}^N \left(X_k -m_k \right) \xi \right) \right] d\xi \right|  \\
& \qquad \leq \tilde{C} \ C(f,g) \left( \frac{ |\supp f | + |\supp g |  }{N^{1- \varepsilon}} + \exp\left(-C\text{dist}\left( \supp f, \supp g \right) \right) \right),
\end{align}
where~$C(f,g)$ is defined by~\eqref{e_def_c}
\end{lemma}
\medskip 

\begin{lemma} [Estimation of the outer integral] \label{l_exp_decay_outer}
Fix~$\delta>0$ small enough. For each~$\varepsilon >0$, there exist constants~$N_2  \in \mathbb{N}$ and~$\tilde{C}= \tilde{C} (\epsilon) \in (0, \infty) $ independent of the external fields~$s=(s_i)$, the boundary~$x^{\mathbb{Z} \setminus [N]}$, and the mean spin~$m$ such that for any~$N \geq N_2$,
\begin{align}
&\left| \int_{\{ \left| \left(1/ \sqrt{N}\right)\xi\right|>\delta\}} \mathbb{E}_{\mu^{\sigma}  } \left[ f(X)g(X) \exp\left( i\frac{1}{\sqrt{N}} \sum_{k=1}^N \left(X_k -m_k \right) \xi \right) \right] d\xi \right| \\
& \qquad \leq \tilde{C} \ C(f,g) \left( \frac{ |\supp f | + |\supp g |  }{N^{1- \varepsilon}} + \exp\left(-C\text{dist}\left( \supp f, \supp g \right) \right) \right),
\end{align}
where~$C(f,g)$ is defined by~\eqref{e_def_c}
\end{lemma}
\medskip

Again, the proof of Proposition~\ref{p_exp_decay_main_computation} is a direct consequence of Lemma~\ref{l_exp_decay_inner} and Lemma~\ref{l_exp_decay_outer}. \\

\noindent \emph{Proof of Proposition~\ref{p_exp_decay_main_computation}.} \ We fix~$\delta>0$ small enough and choose~$N_1, N_2 \in \mathbb{N}$ as in Lemma~\ref{l_exp_decay_inner} and Lemma~\ref{l_exp_decay_outer}. Then for any~$N \geq N_0 = \max(N_1, N_2)$, the triangle inequality followed by Lemma~\ref{l_exp_decay_inner} and Lemma~\ref{l_exp_decay_outer} yields
\begin{align}
&\left| \int_{\mathbb{R}} \mathbb{E}_{\mu ^{\sigma} } \left[ f(X)g(X) \exp\left( i\frac{1}{\sqrt{N}} \sum_{k=1}^N \left(X_k -m_k \right) \xi \right) \right] d\xi \right| \\
& \qquad \leq \left| \int_{\{ \left| \left(1/ \sqrt{N}\right)\xi\right|\leq\delta\}} \mathbb{E}_{\mu^{\sigma}  } \left[ f(X)g(X) \exp\left( i\frac{1}{\sqrt{N}} \sum_{k=1}^N \left(X_k -m_k \right) \xi \right) \right] d\xi \right| \\
& \qquad \quad + \left| \int_{\{ \left| \left(1/ \sqrt{N}\right)\xi\right|>\delta\}} \mathbb{E}_{\mu^{\sigma}  } \left[ f(X)g(X) \exp\left( i\frac{1}{\sqrt{N}} \sum_{k=1}^N \left(X_k -m_k \right) \xi \right) \right] d\xi \right| \\
& \qquad \leq \tilde{C} \ C(f,g) \left( \frac{ |\supp f | + |\supp g |  }{N^{1- \varepsilon}} + \exp\left(-C\text{dist}\left( \supp f, \supp g \right) \right) \right)
\end{align}
\qed

\medskip 

In this section, we shall only provide the proof of Lemma~\ref{l_exp_decay_inner}. The argument presented in the proof of Lemma~\ref{l_equiv_observables_outer} also applies to Lemma~\ref{l_exp_decay_outer}. As in the previous section, we begin with defining auxiliary sets and introduce an auxiliary statement that will be needed in the proof of Lemma~\ref{l_exp_decay_outer}. \\

\begin{figure}[t]
\centering
\begin{tikzpicture}[xscale=1.3]

\draw[fill] (0,0) circle [radius=0.05];
\draw[fill] (.25,0) circle [radius=0.05];
\draw[fill] (.5,0) circle [radius=0.05];

\draw[fill] (0.975,0) circle [radius=0.02];
\draw[fill] (1.125,0) circle [radius=0.02];
\draw[fill] (1.275,0) circle [radius=0.02];

\draw[fill] (1.75,0) circle [radius=0.05];
\draw[fill] (2,0) circle [radius=0.05];
\draw[fill] (2.25,0) circle [radius=0.05];

\draw (2.5,0) circle[radius=0.05];
\draw (2.75,0) circle[radius=0.05];
\draw (3,0) circle[radius=0.05];
\draw (3.25,0) circle[radius=0.05];
\draw (3.5,0) circle[radius=0.05];

\draw[fill] (3.75,0) circle [radius=0.05];
\draw[fill] (4,0) circle [radius=0.05];
\draw[fill] (4.25,0) circle [radius=0.05];
\draw[fill] (4.5,0) circle [radius=0.05];
\draw[fill] (4.75,0) circle [radius=0.05];
\draw[fill] (5,0) circle [radius=0.05];
\draw[fill] (5.25,0) circle [radius=0.05];

\draw (5.5,0) circle[radius=0.05];
\draw (5.75,0) circle[radius=0.05];
\draw (6,0) circle[radius=0.05];
\draw (6.25,0) circle[radius=0.05];
\draw (6.5,0) circle[radius=0.05];

\draw[fill] (6.75,0) circle [radius=0.05];
\draw[fill] (7,0) circle [radius=0.05];
\draw[fill] (7.25,0) circle [radius=0.05];

\draw[fill] (7.725,0) circle [radius=0.02];
\draw[fill] (7.875,0) circle [radius=0.02];
\draw[fill] (8.025,0) circle [radius=0.02];

\draw[fill] (8.5,0) circle [radius=0.05];
\draw[fill] (8.75,0) circle [radius=0.05];
\draw[fill] (9,0) circle [radius=0.05];

\node[align=center, below] at (0,-.05) {$1$};
\node[align=center, below] at (3,-.05) {$i$};
\node[align=center, below] at (6,-.05) {$j$};
\node[align=center, below] at (9,-.05) {$N$};

\draw[decorate,decoration={brace,mirror},thick] (0,-.5) -- node[below]{$F_{2}^{f,g}$} (2.25,-.5);
\draw[decorate,decoration={brace,mirror},thick] (2.5,-.5) -- node[below]{$F_{1}^{f,g}$} (3.5,-.5);
\draw[decorate,decoration={brace,mirror},thick] (3.75,-.5) -- node[below]{$F_{2}^{f,g}$} (5.25,-.5);
\draw[decorate,decoration={brace,mirror},thick] (5.5,-.5) -- node[below]{$F_{1}^{f,g}$} (6.5,-.5);
\draw[decorate,decoration={brace,mirror},thick] (6.75,-.5) -- node[below]{$F_{2}^{f,g}$} (9,-.5);

\draw[thick,<->] (2.5,.5) -- (2.95,.5);
\node[align=center, above] at (2.75,.5) {$L$};
\draw[thick,<->] (3.05,.5) -- (3.5,.5);
\node[align=center, above] at (3.25,.5) {$L$};
\draw[thick,<->] (5.5,.5) -- (5.95,.5);
\node[align=center, above] at (5.75,.5) {$L$};
\draw[thick,<->] (6.05,.5) -- (6.5,.5);
\node[align=center, above] at (6.25,.5) {$L$};
\end{tikzpicture}
\caption{The sets $F_1^{f,g}$ and $F_2 ^{f,g}$ where~$f(x)=x_i-m_i$ and~$g(x)=x_j-m_j$.}\label{f_sets_F}
\end{figure}
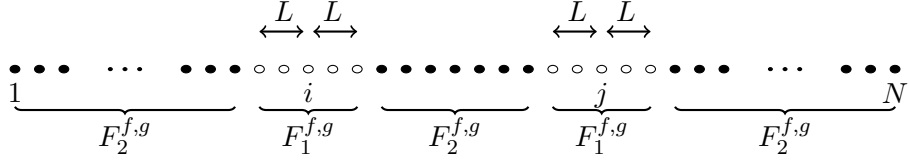

We define the auxiliary sets~$F_1^{f,g}$ and~$F_2^{f,g}$ as (see Figure~\ref{f_sets_F}):
\begin{align}
&F_1^{f,g} : = \{1, 2, \cdots, N \} \cap \{ k : \ \text{dist}\left( k, \supp(f)\right) \leq L \text{ or } \text{dist}\left( k, \supp(g)\right) \leq L   \}, \\
&F_2^{f,g} : = \{1, 2, \cdots, N \} \cap \{ k : \ \text{dist}\left( k, \supp(f)\right) > L \text{ and } \text{dist}\left( k, \supp(g)\right) > L   \},
\end{align}
where~$L\ll N$ is a positive integer that will be chosen later. The main ingredients for this part are Theorem~\ref{p_decay_of_correlations_gce}, Lemma~\ref{l_amgm}, Lemma~\ref{l_fourth_moment} and the extension of~\cite[Lemma 7]{KwMe18}.

\begin{lemma}[Extension of Lemma 7 in~\cite{KwMe18}] \label{l_exchanging_exponential_terms}
For large enough~$N$ and~$\delta>0$ small enough, there exists a positive constant~$C>0$ such that the following inequalities hold for all~$\xi \in \mathbb{R}$ with~$ \frac{\left|\xi\right|}{\sqrt{N}}  \leq \delta$.
\begin{align}
 \left|\mathbb{E}_{\mu^{\sigma}}\left[\exp\left( i \sum_{k \in F_2^{f,g} } \left( X_k -m_k \right) \frac{\xi}{\sqrt{N}}\right) \mid \mathscr{F}_{f,g} \right] \right| \leq C \left( 1 + \xi^2 \right)\exp\left(-C\xi^2 \right) .
\label{e_lemma exchanging 2nd der}
\end{align}
where~$\mathscr{F}_{f,g}$ denotes the sigma algebra defined by
\begin{align}
\mathscr{F}_{f,g} : = \sigma\left(X_k , \ k \in F_1 ^{f,g} \right).
\end{align} 
\end{lemma}
\medskip

For the proof of Lemma~\ref{l_exchanging_exponential_terms} we refer to the proof of~\cite[Lemma 7]{KwMe18} (see also Remark~\ref{r_exchanging_exponential_terms}).  Now we are ready to provide the proof of Lemma~\ref{l_exp_decay_inner}. \\


\noindent \emph{Proof of Lemma~\ref{l_exp_decay_inner}.} \ We define~$e: \mathbb{R}^2  \to \mathbb{C}$ by
\begin{align}
e\left(\xi_1, \xi_2 \right) 
&:= \mathbb{E}_{\mu^{\sigma} } \left[ f(X)g(X) \exp\left( i \sum_{k \in F_1^{f, g} } \left(X_k -m_k \right) \xi_1  +  i \sum_{l \in F_2^{f, g} } \left(X_l -m_l \right) \xi_2\right) \right].
\end{align}
Then a second order Taylor expansion with respect to the first variable~$\xi_1$ yields
\newpage 

\begin{align}
e\left(\xi_1, \xi_2 \right) &= \mathbb{E}_{\mu ^{\sigma}} \left[ f(X)g(X) \exp\left(   i \sum_{l \in F_2^{f, g} } \left(X_l -m_l \right) \xi_2\right) \right] \\
& \quad + \mathbb{E}_{\mu^{\sigma}} \left[ f(X)g(X) \left( i  \sum_{k \in F_1 ^{f,g}} \left(X_k -m_k \right) \right)  \exp\left(  i \sum_{l \in F_2^{f, g} } \left(X_l -m_l \right) \xi_2\right) \right] \xi_1\\
& \quad + \mathbb{E}_{\mu^{\sigma}} \left[ f(X)g(X) \left( i  \sum_{k_1 \in F_1 ^{f,g}} \left(X_{k_1} -m_{k_1} \right) \right)^2 \right. \\
& \qquad \quad \quad  \quad  \left. \times  \exp\left(i \sum_{k_2 \in F_1^{f, g} } \left(X_{k_2} -m_{k_2} \right) \tilde{\xi}_1+   i \sum_{l \in F_2^{f, g} } \left(X_l -m_l \right) \xi_2\right) \right] \xi_1 ^2 ,
\end{align}
where~$\tilde{\xi}_1$ is a real number between~$0$ and~$\xi_1$. In particular for~$\left(\xi_1, \xi_2 \right) = \left( \frac{\xi}{\sqrt{N}}, \frac{\xi}{\sqrt{N}} \right)$, it holds that
\begin{align}
&\mathbb{E}_{\mu^{\sigma} } \left[ f(X)g(X) \exp\left( i \sum_{k=1}^N \left(X_k -m_k \right) \frac{\xi}{\sqrt{N}} \right) \right] \label{e_second_der_a} \\
& \qquad = \mathbb{E}_{\mu^{\sigma} } \left[ f(X)g(X) \exp\left(   i   \sum_{l \in F_2^{f, g} } \left(X_l -m_l \right) \frac{\xi}{\sqrt{N}} \right) \right] \label{e_second_der_a1} \\
& \qquad \quad + \mathbb{E}_{\mu^{\sigma} } \left[ f(X)g(X) \left( i   \sum_{k \in F_1 ^{f,g}} \left(X_k -m_k \right) \right)   \exp\left(  i \sum_{l \in F_2^{f, g} } \left(X_l -m_l \right) \frac{\xi}{\sqrt{N}}\right) \right] \frac{\xi}{\sqrt{N}} \label{e_second_der_a2} \\
& \qquad \quad + \mathbb{E}_{\mu^{\sigma} } \left[ f(X)g(X) \left( i   \sum_{k_1 \in F_1 ^{f,g}} \left(X_{k_1} -m_{k_1} \right) \right)^2 \right.\\
& \qquad \qquad  \quad \quad  \left. \times  \exp\left(i \sum_{k_2 \in F_1^{f, g} } \left(X_{k_2} -m_{k_2} \right) \frac{\tilde{\xi}}{\sqrt{N}}+   i \sum_{l \in F_2^{f, g} } \left(X_l -m_l \right) \frac{\xi}{\sqrt{N}}\right) \right] \left(\frac{\xi}{\sqrt{N}}\right)^2 \label{e_second_der_a3},
\end{align}
where~$\frac{\tilde{\xi}}{\sqrt{N}}$ is a real number between~$0$ and~$\frac{\xi}{\sqrt{N}}$.  

\newpage 

\noindent \textbf{Step 1.} Estimation of~\eqref{e_second_der_a1}. \\ 

By definition of covariances, it holds that
\begin{align}
T_{\eqref{e_second_der_a1}} &= \cov_{\mu^{\sigma}  } \left( f(X)g(X) , \exp\left( i  \sum_{l \in F_2 ^{f,g}} \left(X_l -m_l \right) \frac{\xi}{\sqrt{N}} \right) \right) \label{e_second_der_a11} \\
& \quad + \mathbb{E}_{\mu^{\sigma} } \left[ f(X)g(X) \right] \mathbb{E}_{\mu^{\sigma} } \left[ \exp\left( i  \sum_{l \in F_2 ^{f, g}} \left(X_l -m_l \right) \frac{\xi}{\sqrt{K}} \right)  \right]. \label{e_second_der_a12}
\end{align}
Then a combination of Theorem~\ref{p_decay_of_correlations_gce} and Lemma~\ref{l_amgm} yields
\begin{align}
\left| T_{\eqref{e_second_der_a11}} \right| &\lesssim \|\nabla \left( fg \right) \|_{L^2(\mu^{\sigma})}  \left|\xi\right| \exp\left(-CL\right).
\end{align}
It also follows from Theorem~\ref{p_decay_of_correlations_gce}, Lemma~\ref{l_amgm} and Lemma~\ref{l_exchanging_exponential_terms} that
\begin{align}
\left| T_{\eqref{e_second_der_a12}} \right| & \leq \left| \cov_{\mu^{\sigma} } \left( f(X) , g(X) \right)\right| \left(1+\xi^2 \right) \exp\left(-C\xi^2\right) \\
& \lesssim \|\nabla f \|_{L^2(\mu^{\sigma})}\|\nabla g \|_{L^2(\mu^{\sigma})}   \exp\left(-C \text{dist}\left( \supp f, \supp g \right) \right) \left(1+\xi^2 \right) \exp\left(-C\xi^2 \right).
\end{align}
Therefore recalling the definition~\eqref{e_def_c} of~$C(f,g)$, we obtain
\begin{align}
\left| T_{\eqref{e_second_der_a1}} \right| &\leq \left| T_{\eqref{e_second_der_a11}} \right| + \left| T_{\eqref{e_second_der_a12}} \right| \\
& \lesssim \| \nabla \left( fg \right) \|_{L^2(\mu^{\sigma})} \left(1+ \xi^2 \right) \exp\left(-CL \right) \\
& \quad + \|\nabla f \|_{L^2(\mu^{\sigma})}\|\nabla g \|_{L^2(\mu^{\sigma})}   \exp\left(-C \text{dist}\left( \supp f, \supp g \right) \right) \left(1+\xi^2 \right) \exp\left(-C\xi^2 \right)  \\
& \overset{Holder}{\lesssim} C(f,g)\left(1+ \xi^2 \right) \exp\left(-CL \right) \\
& \quad + C(f,g) \exp\left(-C \text{dist}\left( \supp f, \supp g \right) \right) \left(1+\xi^2 \right) \exp\left(-C\xi^2 \right),
\end{align}
where~$C(f,g)$ is the term defined by~\eqref{e_def_c}. \\

\noindent \textbf{Step 2.} Estimation of~\eqref{e_second_der_a2}. \\ 

Let us begin with introducing auxiliary sets. Denote~$\tilde{F}_1^{f,g}$ and~$\tilde{F}_2^{f,g}$ by
\begin{align}
\tilde{F}_1^{f,g} : = \{ k \in F_1 ^{f, g} : \ \text{dist}\left( k, \supp f \right) \leq \frac{L}{2} \text{ or } \text{dist}\left( k, \supp g \right) \leq \frac{L}{2} \}, \\
\tilde{F}_2^{f,g} : = \{ k \in F_1 ^{f, g} : \ \text{dist}\left( k, \supp f \right) > \frac{L}{2} \text{ and } \text{dist}\left( k, \supp g \right) > \frac{L}{2} \},
\end{align}
and decompose~$T_{\eqref{e_second_der_a2}}$ as
\begin{align}
T_{\eqref{e_second_der_a2}} &= i \frac{\xi}{\sqrt{N}}  \mathbb{E}_{\mu^{\sigma}} \left[ f(X) g(X) \left(\sum_{k \in \tilde{F}_1^{f,g}}\left(X_k -m_k \right)\right) \exp\left(   i   \sum_{l \in F_2^{f, g} } \left(X_l -m_l \right) \frac{\xi}{\sqrt{N}} \right) \right] \label{e_second_der_a21} \\
&\quad + i \frac{\xi}{\sqrt{N}}   \mathbb{E}_{\mu^{\sigma}} \left[ f(X) g(X) \left(\sum_{k \in \tilde{F}_2^{f,g}}\left(X_k -m_k \right)\right) \exp\left(   i   \sum_{l \in F_2^{f, g} } \left(X_l -m_l \right) \frac{\xi}{\sqrt{N}} \right) \right]. \label{e_second_der_a22}
\end{align}
We estimate the terms~\eqref{e_second_der_a21} and~\eqref{e_second_der_a22} separately. \\

\noindent \textbf{Step 2-1.} Estimation of~\eqref{e_second_der_a21}. \\ 

Using the definition of covariance we further decompose the expectation appearing in~$T_{\eqref{e_second_der_a21}}$ as follows:
\begin{align}
&\mathbb{E}_{\mu^{\sigma} } \left[ f(X)g(X) \left(\sum_{k \in \tilde{F}_1^{f,g}}\left(X_k -m_k \right)\right) \exp\left(  i \sum_{l \in F_2^{f, g} } \left(X_l -m_l \right) \frac{\xi}{\sqrt{N}}\right) \right] \label{e_second_der_a21_summand}\\
& = \cov_{\mu^{\sigma} } \left( f(X)g(X) \left(\sum_{k \in \tilde{F}_1^{f,g}}\left(X_k -m_k \right)\right) \ , \ \exp\left(  i \sum_{l \in F_2^{f, g} } \left(X_l -m_l \right) \frac{\xi}{\sqrt{N}}\right)    \right)  \label{e_second_der_a211}  \\
& \quad + \mathbb{E}_{\mu^{\sigma} } \left[ f(X)g(X)\left(\sum_{k \in \tilde{F}_1^{f,g}}\left(X_k -m_k \right)\right) \right] \mathbb{E}_{\mu^{\sigma} } \left[ \exp\left(  i \sum_{l \in F_2^{f, g} } \left(X_l -m_l \right) \frac{\xi}{\sqrt{N}}\right) \right]. \label{e_second_der_a212} 
\end{align}
Let us begin with the estimation of~\eqref{e_second_der_a211}. An application of Theorem~\ref{p_decay_of_correlations_gce} yields
\begin{align}
\left|T_{\eqref{e_second_der_a211}}\right| &\lesssim \left\|\nabla \left( f(X)g(X)\left(\sum_{k \in \tilde{F}_1^{f,g}}\left(X_k -m_k \right)\right)\right) \right\|_{L^2(\mu^{\sigma})} \left|\xi \right| \exp\left(-CL\right) \\
& \leq  \left\| \nabla \left( f(X)g(X)\right) \left(\sum_{k \in \tilde{F}_1^{f,g}}\left(X_k -m_k \right)\right) \right\|_{L^2(\mu^{\sigma})}\left|\xi \right| \exp\left(-CL\right) \label{e_second_der_a2111} \\
& \quad +|\tilde{F}_1^{f,g}|^{\frac{1}{2}}  \cdot  \| f(X)g(X)\|_{L^2(\mu^{\sigma})}  \left|\xi \right| \exp\left(-CL\right). \label{e_second_der_a2112}
\end{align}
Then a combination of H\"{o}lder's inequality and Lemma~\ref{l_fourth_moment} implies that
\begin{align}
T_{\eqref{e_second_der_a2111}} &\lesssim  \| \nabla \left( fg \right)\|_{L^{4}(\mu^{\sigma})} \mathbb{E}_{\mu^{\sigma}} \left[\left(\sum_{k \in \tilde{F}_1^{f,g}}\left(X_k -m_k \right) \right)^4 \right]^{\frac{1}{4}}\left|\xi \right| \exp\left(-CL\right) \\
& \lesssim  C(f,g) |\tilde{F}_1 ^{f,g} |^{\frac{1}{2}} \left|\xi \right| \exp\left(-CL\right).
\end{align}
Next, we use Poincar\'e inequality to obtain
\begin{align}
 \| f(X)g(X)\|_{L^2(\mu^{\sigma})} ^2  &= \var_{\mu^{\sigma}} \left( f(X) g(X) \right) + \left( \mathbb{E}_{\mu^{\sigma}} \left[ f(X)g(X) \right] \right)^2 \\
 & \lesssim \| \nabla \left( f(X) g(X) \right) \|_{L^2 (\mu^{\sigma})} ^2 + \left( \| f(X) \|_{L^2 (\mu^{\sigma})} \| g(X) \|_{L^2 (\mu^{\sigma})}\right)^2 \\
 & \lesssim \| \nabla \left( f(X) g(X) \right) \|_{L^2 (\mu^{\sigma})} ^2  + \| \nabla f \|_{L^2 (\mu^{\sigma})} ^2 \| \nabla g \|_{L^2 (\mu^{\sigma})} ^2. \label{e_l2_fg_estimate}
\end{align}
Therefore we have
\begin{align}
T_{\eqref{e_second_der_a2112}} \lesssim C(f,g) | \tilde{F}_1 ^{f,g} |^{\frac{1}{2}} |\xi | \exp \left(-CL \right).
\end{align}

\medskip 

To estimate the term~\eqref{e_second_der_a212}, we note by Pigeon hole principle that
\begin{align}
\max \left(  \text{dist} \left( k , \supp f \right), \text{dist} \left( k , \supp g \right) \right) \geq \frac{ \text{dist} \left( \supp f , \supp g \right)}{2}.
\end{align}
Let~$S_f$ and~$S_g$ be the sets such that
\begin{align}
S_f : &= \left\{ k \in \tilde{F}_1^{f,g} : \ \text{dist} \left( k , \supp f \right) \geq \frac{ \text{dist} \left( \supp f , \supp g \right)}{2} \right\}, \\
S_g : &= \tilde{F}_1^{f,g} \setminus S_f.
\end{align}
Further decomposing the term~$T_{\eqref{e_second_der_a212}}$, we have
\begin{align}
\left| T_{\eqref{e_second_der_a212}} \right| &\overset{Lemma~\ref{l_exchanging_exponential_terms}}{\lesssim} \left| \cov_{\mu^{\sigma}} \left( f(X) \ , \ g(X)\left(\sum_{k \in S_f}\left(X_k -m_k \right)\right) \right) \right| \left(1+ \xi^2 \right) \exp\left(-C\xi^2 \right) \\
& \qquad \quad +  \left| \cov_{\mu^{\sigma}} \left( f(X)\left(\sum_{k \in S_g}\left(X_k -m_k \right)\right) \ , \ g(X) \right) \right| \left(1+ \xi^2 \right) \exp\left(-C\xi^2 \right) \\
& \overset{Theorem~\ref{p_decay_of_correlations_gce}}{\lesssim} \| \nabla f \|_{L^2 (\mu^{\sigma})} \left\| \nabla \left(  g(X) \left( \sum_{k \in S_f} \left( X_k -m_k \right) \right) \right) \right\|_{L^2 (\mu^{\sigma})} \\
& \qquad \qquad \times \exp \left(- C \text{dist} ( \supp f \ , \ \supp g) \right) (1 + \xi^2 ) \exp \left( -C \xi ^2 \right)  \label{e_a212_sf}\\
& \qquad \quad + \| \nabla g \|_{L^2 (\mu^{\sigma})} \left\| \nabla \left(  f(X) \left( \sum_{k \in S_g} \left( X_k -m_k \right) \right) \right) \right\|_{L^2 (\mu^{\sigma})} \\
& \qquad \qquad \quad  \times \exp \left(- C \text{dist} ( \supp f \ , \ \supp g) \right) (1 + \xi^2 ) \exp \left( -C \xi ^2 \right). \label{e_a212_sg}
\end{align}
Note that we have
\begin{align}
&\left\| \nabla \left(  g(X) \left( \sum_{k \in S_f} \left( X_k -m_k \right) \right) \right) \right\|_{L^2 (\mu^{\sigma})} \\
&\qquad  \leq \left\| \nabla g(X)\left( \sum_{k \in S_f} \left( X_k -m_k \right) \right) \right\|_{L^2(\mu^{\sigma})} + |S_f |^{\frac{1}{2}} \cdot \| g \|_{L^2 (\mu^{\sigma})} \\
& \quad \quad \overset{Holder}{\leq} \| \nabla g \|_{L^{4}(\mu^{\sigma})} \left\| \sum_{k \in S_f} \left( X_k -m_k\right) \right\|_{L^{4} (\mu^{\sigma})} + |S_f |^{\frac{1}{2}}\|g \|_{L^2 (\mu^{\sigma})} \\
& \overset{Lemma~\ref{l_fourth_moment}, Poincare , Holder}{\lesssim} |S_f |^{\frac{1}{2}} \|\nabla g \|_{L^{4}(\mu^{\sigma})}. \label{e_further_estimation2}
\end{align}
Using similar estimation as in~\eqref{e_further_estimation2}, it holds that
\begin{align}
\left| T_{\eqref{e_second_der_a212}} \right| &\leq \left| T_{\eqref{e_a212_sf}} \right|  +\left| T_{\eqref{e_a212_sg}} \right|   \\
& \lesssim C(f,g) \left( |S_f |^{\frac{1}{2}} + |S_g|^{\frac{1}{2}} \right)\exp \left(- C \text{dist} ( \supp f \ , \ \supp g) \right) (1 + \xi^2 ) \exp \left( -C \xi ^2 \right) \\
& \lesssim C(f,g) |\tilde{F}_1 ^{f,g}|^{\frac{1}{2}} \exp \left(- C \text{dist} ( \supp f \ , \ \supp g) \right) (1 + \xi^2 ) \exp \left( -C \xi ^2 \right).
\end{align}
As a consequence we obtain
\begin{align}
\left| T_{\eqref{e_second_der_a21_summand}}\right| &\leq \left| T_{\eqref{e_second_der_a211}} \right|+\left| T_{\eqref{e_second_der_a212}} \right|  \\
&\lesssim C(f,g) |\tilde{F}_1 ^{f,g} |^{\frac{1}{2}} \left|\xi \right| \exp\left(-CL\right) \\
&\quad + C(f,g) |\tilde{F}_1 ^{f,g}|^{\frac{1}{2}} \exp \left(- C \text{dist} ( \supp f \ , \ \supp g) \right) (1 + \xi^2 ) \exp \left( -C \xi ^2 \right).
\end{align}
To conclude, we plug in this estimate into~\eqref{e_second_der_a21} and get
\begin{align}
\left| T_{\eqref{e_second_der_a21}}\right| &\lesssim C(f,g) \xi^2  \exp \left(-CL\right) \frac{| \tilde{F}_1^{f,g}|^{\frac{1}{2}}}{\sqrt{N}} \\
& \quad + C(f,g)\exp\left(-C \text{dist}\left(\supp(f), \supp(g) \right) \right) \left|\xi\right|\left(1+ \xi^2 \right) \exp\left(-C\xi^2 \right)\frac{| \tilde{F}_1^{f,g}|^{\frac{1}{2}}}{\sqrt{N}} .
\end{align}
\medskip 

\noindent \textbf{Step 2-2.} Estimation of~\eqref{e_second_der_a22}. \\

We write the expectation term in~$T_{\eqref{e_second_der_a22}}$ as
\begin{align}
&\mathbb{E}_{\mu^{\sigma}} \left[ f(X)g(X)  \left(\sum_{k \in \tilde{F}_2 ^{f,g}} \left(X_k -m_k \right)\right) \exp\left(  i \sum_{l \in F_2^{f, g} } \left(X_l -m_l \right) \frac{\xi}{\sqrt{N}}\right) \right] \label{e_second_der_a22_summand} \\
& \qquad = \cov_{\mu^{\sigma}} \left( f(X)g(X) , \left(\sum_{k \in \tilde{F}_2 ^{f,g}} \left(X_k -m_k \right)\right) \exp\left(  i \sum_{l \in F_2^{f, g} } \left(X_l -m_l \right) \frac{\xi}{\sqrt{N}}\right)    \right)  \label{e_second_der_a221} \\
& \qquad \quad + \mathbb{E}_{\mu^{\sigma}} \left[ f(X)g(X)  \right] \mathbb{E}_{\mu^{\sigma} } \left[ \left(\sum_{k \in \tilde{F}_2 ^{f,g}} \left(X_k -m_k \right)\right)\exp\left(  i \sum_{l \in F_2^{f, g} } \left(X_l -m_l \right) \frac{\xi}{\sqrt{N}}\right) \right]. \label{e_second_der_a222}
\end{align}
Then applying similar arguments as in~\textbf{Step 2-1} gives
\begin{align}
&\left| T_{\eqref{e_second_der_a221} } \right| \lesssim C(f,g) | \tilde{F}_2 ^{f,g }| ^{\frac{1}{2}} \left(1+ \xi^2 \right) \exp\left(-CL \right), \\
&\left| T_{\eqref{e_second_der_a222} } \right| \lesssim C(f,g)| \tilde{F}_2 ^{f,g }| ^{\frac{1}{2}} \exp\left(-C \text{dist}\left( \supp f, \supp g \right) \right) \left(1+\xi^2 \right) \exp\left(-C\xi^2 \right).
\end{align}
Therefore we have
\begin{align}
\left|T_{\eqref{e_second_der_a22_summand}}\right| &\leq \left| T_{\eqref{e_second_der_a221} } \right| + \left| T_{\eqref{e_second_der_a222} } \right| \\
& \lesssim C(f,g)| \tilde{F}_2 ^{f,g }| ^{\frac{1}{2}} \left(1 + \xi^2 \right) \exp \left(-CL\right) \\
&\quad + C(f,g)| \tilde{F}_2 ^{f,g }| ^{\frac{1}{2}} \exp\left(-C \text{dist}\left(\supp(f), \supp(g) \right) \right) \left(1+ \xi^2 \right) \exp\left(-C\xi^2 \right),
\end{align}
and as a consequence
\begin{align}
\left|T_{\eqref{e_second_der_a22}} \right| &\lesssim  C(f,g) \left| \xi \right| \left(1 + \xi^2 \right) \exp \left(-CL\right) \frac{| \tilde{F}_2^{f,g}|^{\frac{1}{2}}}{\sqrt{N}} \\
& \quad + C(f,g)\exp\left(-C \text{dist}\left(\supp(f), \supp(g) \right) \right) \left|\xi\right|\left(1+ \xi^2 \right) \exp\left(-C\xi^2 \right)\frac{| \tilde{F}_2^{f,g}|^{\frac{1}{2}}}{\sqrt{N}} .
\end{align}

\medskip 

\noindent \textbf{Step 2-3.} Conclusion of \textbf{Step 2}. \\

From the two steps (\textbf{Step 2-1} and~\textbf{Step 2-2}) we obtain enough estimate
\begin{align}
\left|T_{\eqref{e_second_der_a2}}\right| &\leq \left| T_{\eqref{e_second_der_a21}} \right| +  \left| T_{\eqref{e_second_der_a22}} \right| \\
&\lesssim C(f,g) \left| \xi \right| \left(1 + \xi^2 \right) \exp \left(-CL\right) \frac{| F_1^{f,g}|^{\frac{1}{2}}}{\sqrt{N}} \\
& \quad + C(f,g)\exp\left(-C \text{dist}\left(\supp(f), \supp(g) \right) \right) \left|\xi\right|\left(1+ \xi^2 \right) \exp\left(-C\xi^2 \right)\frac{| F_1^{f,g}|^{\frac{1}{2}}}{\sqrt{N}},
\end{align}
where we used the trivial inequality
\begin{align}
|\tilde{F}_1 ^{f,g} | , |\tilde{F}_2 ^{f,g} | \leq |F_1 ^{f,g}|.
\end{align}
\medskip

\noindent \textbf{Step 3.} Estimation of~\eqref{e_second_der_a3}. \\

We consider the conditional expectation with respect to the sigma algebra~$\mathscr{F}_{f,g} = \sigma \left( X_k , \  k \in F_1^{f,g} \right)$. Then it follows from Lemma~\ref{l_exchanging_exponential_terms} that 
\begin{align}
T_{\eqref{e_second_der_a3}} &= \left|\mathbb{E}_{\mu^{\sigma}} \left[ f(X)g(X) \left( i \sum_{k_1 \in F_1 ^{f, g} } \left(X_{k_1} -m_{k_1} \right) \right)^2   \exp\left( i \sum_{k_2 \in F_1 ^{f,g}} \left(X_{k_2} -m_{k_2} \right) \frac{\xi}{\sqrt{N}} \right) \right. \right. \\
& \left. \left. \qquad \qquad \times \mathbb{E}_{\mu^{\sigma}} \left[ \exp\left(  \sum_{l \in F_2 ^{f,g}} \left(X_l -m_l \right) \frac{\xi}{\sqrt{N}} \right) \mid \mathscr{F}_{f,g}  \right] \right]  \right| \frac{\xi^2}{N} \\
& \overset{Lemma~\ref{l_exchanging_exponential_terms}}{\lesssim} \mathbb{E}_{\mu^{\sigma}} \left| f(X)g(X) \left( i \sum_{k_1 \in F_1 ^{f, g} } \left(X_{k_1} -m_{k_1} \right) \right) ^2  \right|  \left(1+ \xi^2 \right)\exp\left(-C\xi^2 \right) \frac{\xi^2}{N} \\
& \overset{Holder}{\lesssim} \|f(X) g(X) \|_{L^2 ( \mu^{\sigma})} \mathbb{E}_{\mu^{\sigma}} \left[ \left( \sum_{k_1 \in F_1 ^{f,g}} \left( X_{k_1} - m_{k_1} \right)\right)^4\right]^{\frac{1}{2}}   \left(1+ \xi^2 \right) \exp\left(-C\xi^2 \right)\frac{\xi^2}{N}.
\end{align}
Now we use~\eqref{e_l2_fg_estimate} and apply Lemma~\ref{l_fourth_moment} to get
\begin{align}
T_{\eqref{e_second_der_a3}} \lesssim C(f,g) \frac{| F_1 ^{f,g} |}{N} \xi^2 (1 + \xi^2) \exp \left( -C \xi^2 \right) .
\end{align}

\medskip 

\noindent \textbf{Step 4.} Conclusion. \\

To conclude, we have
\begin{align}
\left| T_{\eqref{e_second_der_a}} \right| & \leq  \left|T_{\eqref{e_second_der_a1}} \right| + \left|T_{\eqref{e_second_der_a2}} \right|+ \left|T_{\eqref{e_second_der_a3}} \right| \\
&\leq  C(f,g)\left(1+ \xi^2 \right) \exp\left(-CL \right) \\
& \quad + C(f,g) \exp\left(-C \text{dist}\left( \supp f, \supp g \right) \right) \left(1+\xi^2 \right) \exp\left(-C\xi^2 \right) \\
& \quad + C(f,g) \left| \xi \right| \left(1 + \xi^2 \right) \exp \left(-CL\right) \frac{| F_1^{f,g}|^{\frac{1}{2}}}{\sqrt{N}} \\
& \quad + C(f,g)\exp\left(-C \text{dist}\left(\supp(f), \supp(g) \right) \right) \left|\xi\right|\left(1+ \xi^2 \right) \exp\left(-C\xi^2 \right)\frac{| F_1^{f,g}|^{\frac{1}{2}}}{\sqrt{N}} \\
& \quad + C(f,g) \frac{| F_1 ^{f,g} |}{N} \xi^2 (1 + \xi^2) \exp \left( -C \xi^2 \right).
\end{align}
Note that for~$L=N^{\varepsilon} \ll N$ and~$N$ large enough, it holds that
\begin{align}
\int_{\{ \left| \left(1/ \sqrt{N}\right)\xi\right|\leq\delta\}} |\xi|^k (1+ \xi^2) \exp\left(-CL\right) d\xi \lesssim \frac{1}{N^2} \qquad \text{for } k=0, 1,\\
\int_{\{ \left| \left(1/ \sqrt{N}\right)\xi\right|\leq\delta\}} |\xi|^k (1+ \xi^2) \exp\left(-C\xi^2\right) d\xi \lesssim 1 \qquad \text{for } k=0, 1, 2.
\end{align}
Trivially, we have~$|F_1^{f,g}| \leq N $ and also
\begin{align}
|F_1 ^{f,g} | &\leq  2L \left( |\supp f | + |\supp g| \right) .
\end{align}
Then a combination of the estimates from above yields the desired estimate
\begin{align}
\left| T_{\eqref{e_exp_decay_inner}} \right| & \leq \tilde{C} \ C(f,g) \left( \frac{ |\supp f | + |\supp g |  }{N^{1- \varepsilon}} + \exp\left(-C\text{dist}\left( \supp f, \supp g \right) \right) \right) . 
\end{align}
where~$C(f,g)$ is defined by~\eqref{e_def_c} 
\qed
\medskip 

\begin{remark}
One could hope that the bound on Proposition~\ref{p_equiv_observables_main_computation} can be improved via second order Taylor expansion as in Proposition~\ref{p_exp_decay_main_computation}. However, this is very hard because the corresponding first order term
\begin{align}
\mathbb{E}_{\mu^{\sigma} } \left[ f(X) \left( i   \sum_{k \in E_1 ^{f}} \left(X_k -m_k \right) \right)   \exp\left(  i \sum_{l \in E_2^{f} } \left(X_l -m_l \right) \frac{\xi}{\sqrt{N}}\right) \right] \frac{\xi}{\sqrt{N}}
\end{align}
does not have enough decay. For example, considering the case~$f(x) = x-m_i$ and $k = i \in E_1 ^f$, we have
\begin{align}
\mathbb{E}_{\mu^{\sigma} } \left[ \left(X_i -m_i \right)^2   \exp\left(  i \sum_{l \in E_2^{f} } \left(X_l -m_l \right)  \frac{\xi}{\sqrt{N}}\right) \right] i \frac{\xi}{\sqrt{N}} \sim \frac{\xi}{\sqrt{N}},
\end{align}
which is insufficient to deduce~$\frac{1}{N^{1-\varepsilon}}$ scale bound.
\end{remark}
\bigskip 

\section{Proof of Theorem~\ref{p_uniqueness_iv_gibbs_ce}}\label{s_uniqueness_ivGm}
\noindent \emph{Proof of Theorem~\ref{p_uniqueness_iv_gibbs_ce}. } \  Suppose that there are two infinite-volume Gibbs measures~$\mu$ and~$\nu$ of the ce~$\mu_m^{\Lambda}$. It suffices to prove that for any smooth function~$f : \mathbb{R}^{\mathbb{Z}} \to \mathbb{R}$ with bounded support
\begin{align}
 \int f \mu = \int f \nu.
\end{align}
Let us fix a smooth function~$f$ with bounded support. For each~$r>0$ define~$B_r \subset \mathbb{Z}$ by
\begin{align}
B_r : = \left\{k \in \mathbb{Z} \ | \ -r < k < r \right\}
\end{align}
and choose~$K>0$ so that
\begin{align} \label{e_suppf_subset}
\supp f \subset B_K.
\end{align}
For each~$r>K$ we decompose the measure~$\mu$ into the conditional measure~$\mu \left( dx^{B_r} | y^{\mathbb{Z} \backslash B_r} \right)$ and the marginal measure~$\bar{\mu}\left(dy^{ \mathbb{Z} \backslash B_r} \right)$, i.e. for any test function g
\begin{align}
\int g \mu = \int \int g\left(x^{B_r}, y^{\mathbb{Z} \backslash B_r} \right)\mu \left( dx^{B_r} | y^{\mathbb{Z} \backslash B_r} \right) \bar{\mu}\left(dy^{ \mathbb{Z} \backslash B_r} \right).
\end{align}
Similarly, decompose the measure~$\nu$ into~$\nu \left( dx^{B_r} | y^{\mathbb{Z} \backslash B_r} \right)$ and~$\bar{\nu}\left(dy^{ \mathbb{Z} \backslash B_r} \right)$. Then it holds from (DLR) equations that
\begin{align} \label{e_dlr_compare}
\mu \left( dx^{B_r} | y^{\mathbb{Z} \backslash B_r} \right) = \nu \left( dx^{B_r} | y^{\mathbb{Z} \backslash B_r} \right)= \mu_m^{B_r}(dx^{B_r} | y^{\mathbb{Z} \backslash B_r} ).
\end{align}
For notational convenience we write
\begin{align}
x= x^{B_r} , \qquad y= y^{\mathbb{Z} \backslash B_r} \qquad \text{and} \qquad z=z^{\mathbb{Z} \backslash B_r}. \label{e_convention}
\end{align}
Note that~\eqref{e_dlr_compare} implies 
\begin{align}
\left| \int f \mu - \int f \nu \right| &  = \left| \int \int f \mu \left(dx | y\right) \bar{\mu} \left(dy \right) - \int \int f \nu \left(dx | z \right) \bar{\nu} \left(dz \right)  \right|\\
&  = \left| \int \int \left( \int f \mu_m^{B_r} \left(dx | y\right)  -  \int f \mu_m^{B_r} \left(dx | z\right)   \right)  \bar{\mu} \left(dy \right) \bar{\nu} \left(dz \right) \right| \\
&  \leq  \int \int \left| \int f \mu_m^{B_r} \left(dx | y\right)  -  \int f \mu_m^{B_r} \left(dx | z\right)   \right|  \bar{\mu} \left(dy \right) \bar{\nu} \left(dz \right).  \label{e_uniqueness_claim}
\end{align}
We claim that the right hand side of~$T_{\eqref{e_uniqueness_claim}}$ becomes small when choosing~$r>0$ large enough. More precisely, we have the following estimate.

\begin{lemma}\label{l_core_estimate_uniqueness}
Let~$\varepsilon$ be a fixed positive number. Then it holds that
\begin{align}
&\int \int \left| \int f \mu_m^{B_r} \left(dx | y\right)  -  \int f \mu_m^{B_r} \left(dx | z\right)   \right|  \bar{\mu} \left(dy \right) \bar{\nu} \left(dz \right) \\
&\qquad  \lesssim R^2 \left( |f|_{\infty} + | \nabla f |_{\infty} \right) \left( \frac{|\supp f |}{r^{1-\varepsilon}} +  \exp\left( -C (r-R-K) \right) \right).
\end{align}
\end{lemma}
\medskip 

The statement from above finishes the proof of Theorem~\ref{p_uniqueness_iv_gibbs_ce} by letting~$r \to \infty$ and get
\begin{align}
\left| \int f \mu - \int f \nu \right| =0.
\end{align}
\qed

\medskip

Now let us turn to the proof of Lemma~\ref{l_core_estimate_uniqueness}. \\

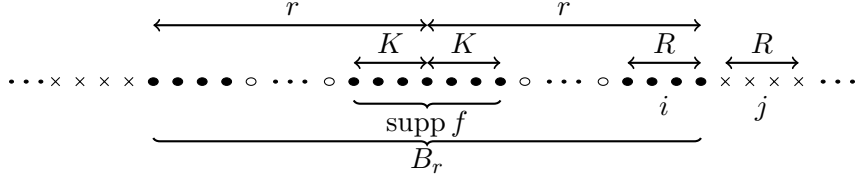
\begin{figure}[t]
\centering
\begin{tikzpicture}[xscale=1.3]
\draw[fill] (0.25,0) circle [radius=0.02];
\draw[fill] (0.4,0) circle [radius=0.02];
\draw[fill] (0.55,0) circle [radius=0.02];
\draw (0.7,0) node[cross] {};
\draw (0.95,0) node[cross] {};
\draw (1.2,0) node[cross] {};
\draw (1.45,0) node[cross] {};

\draw[fill] (1.7,0) circle [radius=0.05];
\draw[fill] (1.95,0) circle [radius=0.05];
\draw[fill] (2.2,0) circle [radius=0.05];
\draw[fill] (2.45,0) circle [radius=0.05];

\draw (2.7,0) circle[radius=0.05];
\draw[fill] (2.95,0) circle [radius=0.02];
\draw[fill] (3.1,0) circle [radius=0.02];
\draw[fill] (3.25,0) circle [radius=0.02];
\draw (3.5,0) circle[radius=0.05];

\draw[fill] (3.75,0) circle [radius=0.05];
\draw[fill] (4,0) circle [radius=0.05];
\draw[fill] (4.25,0) circle [radius=0.05];
\draw[fill] (4.5,0) circle [radius=0.05];
\draw[fill] (4.75,0) circle [radius=0.05];
\draw[fill] (5,0) circle [radius=0.05];
\draw[fill] (5.25,0) circle [radius=0.05];

\draw (5.5,0) circle[radius=0.05];
\draw[fill] (5.75,0) circle [radius=0.02];
\draw[fill] (5.9,0) circle [radius=0.02];
\draw[fill] (6.05,0) circle [radius=0.02];
\draw (6.3,0) circle[radius=0.05];

\draw[fill] (6.55,0) circle [radius=0.05];
\draw[fill] (6.8,0) circle [radius=0.05];
\draw[fill] (7.05,0) circle [radius=0.05];
\draw[fill] (7.3,0) circle [radius=0.05];

\draw (7.55,0) node[cross] {};
\draw (7.8,0) node[cross] {};
\draw (8.05,0) node[cross] {};
\draw (8.3,0) node[cross] {};
\draw[fill] (8.55,0) circle [radius=0.02];
\draw[fill] (8.7,0) circle [radius=0.02];
\draw[fill] (8.85,0) circle [radius=0.02];

\node[align=center, below] at (6.925,-.05) {$i$};
\node[align=center, below] at (7.925,-.05) {$j$};

\draw[decorate,decoration={brace,mirror},thick] (3.75,-.25) -- node[below]{$\supp f$} (5.25,-.25);
\draw[decorate,decoration={brace,mirror},thick] (1.7,-.75) -- node[below]{$B_r$} (7.3,-.75);

\draw[thick,<->] (3.75,.25) -- (4.5,.25);
\node[align=center, above] at (4.125,.25) {$K$};
\draw[thick,<->] (4.5,.25) -- (5.25,.25);
\node[align=center, above] at (4.875,.25) {$K$};
\draw[thick,<->] (6.55,.25) -- (7.3,.25);
\node[align=center, above] at (6.925,.25) {$R$};
\draw[thick,<->] (7.55,.25) -- (8.3,.25);
\node[align=center, above] at (7.925,.25) {$R$};
\draw[thick,<->] (1.7,.75) -- (4.5,.75);
\node[align=center, above] at (3.125,.75) {$r$};
\draw[thick,<->] (4.5,.75) -- (7.3,.75);
\node[align=center, above] at (5.9,.75) {$r$};

\end{tikzpicture}
\caption{Distance between~$\supp f$ and~$\{i\}$}\label{f_supp}
\end{figure}

\noindent \emph{Proof of Lemma~\ref{l_core_estimate_uniqueness}. } \ By interpolation it holds that (recall the convention~\eqref{e_convention})
\begin{align}
\int f \mu_m^{B_r} \left(dx | y\right)  -  \int f \mu_m^{B_r} \left(dx | z\right) &  = \int_0 ^1  \left( \frac{d}{dt} \int f \mu_m ^{B_r} \left(dx | ty+ (1-t)z \right) \right)  dt \\
& = \int_0 ^1 \cov_{\mu_m^{B_r}\left(dx | ty+ (1-t)z  \right)} \left( f, \sum_{ \substack{ i \in B_r, \ j \notin B_r \\ |i-j| \leq R} } M_{ij} x_i (y_j - z_j) \right)dt. \label{e_interpolation}
\end{align}
Let us consider the integrand in~\eqref{e_interpolation}. To estimate the covariance with respect to the measure~$\mu_m^{B_r}\left(dx | ty+ (1-t)z  \right)$, let us define the corresponding gce~$\mu^{B_r, \tau} \left( dx | ty + (1-t)z\right)$ by
\begin{align}
\mu^{B_r, \tau} \left( dx | ty + (1-t)z\right) = \frac{1}{Z} \exp \left( \tau \sum_{k \in B_r} x_k - H\left( x, ty+(1-t)z \right) \right) dx,
\end{align}
where we choose~$\tau = \tau(m)$ such that (cf. Definition~\ref{a_relation_sigma_m})
\begin{align}
m = \frac{1}{|B_r|} \int \left( \sum_{k \in B_r} x_k\right) \mu^{B_r, \tau} \left( dx | ty + (1-t)z\right) .
\end{align}
For a pair~$(i,j)$ with~$i \in B_r$,~$j \notin B_r$ and~$|i-j|\leq R$, the triangle inequality yields
\begin{align}
|i| \geq |j| - |i-j| \geq r - R,
\end{align}
and in particular for~$r> R+ K$ (cf.~\eqref{e_suppf_subset} and Figure~\ref{f_supp}),
\begin{align}
\text{dist} \left( \supp f , \{i \}  \right) \geq r-R-K.
\end{align}

Then a combination of Theorem~\ref{p_exponential_decay_of_correlations_of_canonical_ensemble} and Lemma~\ref{l_amgm} yields
\begin{align}
&\left| \cov_{\mu_m^{B_r}\left(dx | ty+ (1-t)z  \right)} \left( f, \sum_{ \substack{ i \in B_r, \ j \notin B_r \\ |i-j| \leq R} } M_{ij} x_i (y_j - z_j) \right)\right| \\
& \qquad \leq  \sum_{ \substack{ i \in B_r, \ j \notin B_r \\ |i-j| \leq R} } |M_{ij}| |y_j - z_j | \left| \cov_{\mu_m^{B_r}\left(dx | ty+ (1-t)z  \right)} \left( f, x_i \right) \right| \\
& \qquad \lesssim  \left( |f|_{\infty} + | \nabla f |_{\infty} \right) \left( \frac{|\supp f |}{r^{1-\varepsilon}} +  \exp\left( -C (r-R-K) \right) \right) \left( \sum_{ \substack{ i \in B_r, \ j \notin B_r \\ |i-j| \leq R} } |y_j - z_j| \right)  .\label{e_interpolation_estimate}
\end{align}
Hence~\eqref{e_interpolation} and~\eqref{e_interpolation_estimate} imply
\begin{align}
&\int \int \left| \int f \mu_m^{B_r} \left(dx | y\right)  -  \int f \mu_m^{B_r} \left(dx | z\right)   \right|  \bar{\mu} \left(dy \right) \bar{\nu} \left(dz \right)  \\
& \qquad \leq \left( |f|_{\infty} + | \nabla f |_{\infty} \right) \left( \frac{|\supp f |}{r^{1-\varepsilon}} +  \exp\left( -C (r-R-K) \right) \right) \\
& \qquad \quad \times \sum_{ \substack{ i \in B_r, \ j \notin B_r \\ |i-j| \leq R} } \int \int  |y_j - z_j |  \bar{\mu} \left(dy \right) \bar{\nu} \left(dz \right). 
\end{align}
Note that Cauchy's inequality implies
\begin{align}
\int \int  |y_j -z_j | \bar{\mu}(dy) \bar{\mu}(dz) & \leq \left( \int \int \left(y_j-z_j\right) ^2 \bar{\mu}(dy)\bar{\mu}(dz) \right)^{\frac{1}{2}} \\
& = \left( 2 \var_{\bar{\mu}} \left( y_j\right) \right)^{\frac{1}{2}}\overset{\eqref{e_uniform_bound_uniqueness}}{\lesssim} 1.
\end{align}
Because there are at most~$2 R^2$ many pairs of~$(i,j)$ with~$i \in B_r$,~$j \notin B_r$ and~$|i-j| \leq R$, we have
\begin{align}
&\int \int \left| \int f \mu_m^{B_r} \left(dx | y\right)  -  \int f \mu_m^{B_r} \left(dx | z\right)   \right|  \bar{\mu} \left(dy \right) \bar{\nu} \left(dz \right) \\
& \qquad \lesssim R^2 \left( |f|_{\infty} + | \nabla f |_{\infty} \right) \left( \frac{|\supp f |}{r^{1-\varepsilon}} +  \exp\left( -C (r-R-K) \right) \right).
\end{align}
\qed

\section*{Acknowledgment}
This research has been partially supported by NSF grant DMS-1407558. The authors are thankful to many people discussing the problem and helping to improve the preprint. Among them are Tim Austin, Frank Barthe, Marek Biskup, Pietro Caputo, Jean-Dominique Deuschel, Max Fathi, Andrew Krieger, Michel Ledoux, Sangchul Lee, Thomas Liggett, Guido Mont\'ufar, Felix Otto, Andr\'e Schlichting, Daniel Ueltschi, and Tianqi Wu. The authors want to thank Marek Biskup, UCLA and KFAS for financial support.

\bibliographystyle{alpha}
\bibliography{bib}

\end{document}